\documentclass[10pt]{article}
 \usepackage{epsf}

\usepackage[utf8x]{inputenc}
\usepackage[T1]{fontenc}
\usepackage[english]{babel}

\usepackage{media9}

\usepackage{amsmath, amsfonts, amstext, amssymb} 
\usepackage{pythontex}
 \usepackage{color,colortbl}
 \usepackage{fancyhdr}
 \usepackage{xcolor,pict2e}
 \definecolor{myaqua}{rgb}{0.0,0.5,0.55}
 \definecolor{lightaqua}{rgb}{0.75,0.95,0.95}
\usepackage{ epsfig, lscape}
 \usepackage[colorlinks = true,
            linkcolor = myaqua,
            urlcolor  = blue,
            citecolor = red,
            pdfpagemode=UseOutlines,
            bookmarksnumbered=true,
            pdfpagelabels,
            breaklinks]{hyperref}
\usepackage{chemfig}
\usepackage{tikz}
\usetikzlibrary{automata}
\usetikzlibrary{shapes}
\usetikzlibrary{arrows,positioning,automata,shadows}
\usepackage{lipsum}
\usepackage{graphics} 
\usepackage[all,dvips]{xy}
\xyoption{ps} 

\usepackage{tikzscale}
\usepackage{tikz}
\usetikzlibrary{arrows,backgrounds,calc,shapes,patterns}
\usetikzlibrary{er}
\usepackage{tikz-3dplot}

\usepackage{pgfplots}
\usetikzlibrary{pgfplots.groupplots}

\usepackage{pgfplotstable}
\usepackage{currfile}
\usepackage{lmodern}
\usepackage{enumitem}
\usepackage{filecontents}
\usepackage{pdftricks}
\usetikzlibrary{spy}

\usetikzlibrary{external}
\newcommand{\atm}{\entrymodifiers={++[o][F-]}}

\usepackage{rcs}
\usepackage{pgfplots}

\usepackage{pst-node,pst-plot}

\colorlet{bleuturquoise}{cyan}
\usepackage{siunitx}

\usepackage{filecontents}
\usepackage{grffile}
\usepackage{textcomp}

\newcommand{\N}{\mathbb{N}}
\newcommand{\Z}{\mathbb{Z}}
\newcommand{\R}{\mathbb{R}}
\newcommand{\M}{\mathcal{M}}
\newcommand{\A}{\mathcal{A}}

\newcommand{\C}{\mathfrak{C}}

\newtheorem{theo}{Theorem}
\newtheorem{coro}{Corollary}
\newtheorem{lem}{Lemma}
\newtheorem{prop}{Proposition}

\newtheorem{rem}{Remark}

\newenvironment{preu}{\textbf{Proof }\rm}{\par\hfill$\square$}
\newenvironment{preu 1}{\textbf{Proof of Theorem \ref{t2}.}\rm}{\par\hfill$\square$}
\newenvironment{preu 2}{\textbf{Proof of Theorem \ref{t1}.}\rm}{\par\hfill$\square$}
\newenvironment{preu 3}{\textbf{Proof of Theorem \ref{th3}}\rm}{\par\hfill$\square$}
\newenvironment{preu 4}{\textbf{Proof of Theorem \ref{theo3}}\rm}{\par\hfill$\square$}
\newenvironment{preu 5}{\textbf{Proof of Lemma \ref{lmm3} }\rm}{\par\hfill$\square$}
\newtheorem{defi}{Definition}
\newtheorem{exple}{Example}

\usepackage{datetime}
\usepackage{lineno}
\usepackage{floatrow}
\def\lin#1#2{\textcolor[rgb]{0.6,0.6,0.6}{\vspace*{#1mm} \hrule
   height 3 pt \vspace*{#2mm}}}
\def\bt{\begin{tabular}}
\def\et{\end{tabular}}
\def\and{\mbox{ and }}

\def\P{\mbox{\bf P}}

\def\1{{\bf 1}}

\begin{document}

% \fancyhead[L]{\hspace*{-13mm}
% \bt{l}{\bf Open Journal of *****, 2014, *,**}\\
% Published Online **** 2014 in SciRes.
% \href{http://www.scirp.org/journal/*****}{\color{blue}{\underline{\smash{http://www.scirp.org/journal/****}}}} \\
% \href{http://dx.doi.org/10.4236/****.2014.*****}{\color{blue}{\underline{\smash{http://dx.doi.org/10.4236/****.2014.*****}}}} \\
% \et}
 \fancyhead[R]{\includegraphics{pic1.ps}}

 $\mbox{ }$

 \vskip 12mm
 
 { % \fontfamily{Cambria}\selectfont

% "Title of the Paper"
{\noindent{\Large\bf\color{myaqua}
   Exchange of Intervals and Intrinsic Ergodicity of the $\beta$-shift}} %\mathbb{R}_{+} $\Sigma(H)$
%
% \runtitle{Change-Point Analysis of Survival Data}
\\[6mm]
{\bf Florent NGUEMA NDONG$^1$, Anne BERTRAND $^2$}}
%{\bf Florent NGUEMA NDONG$^1$, Jules TINDZOGHO NTSIRI$^2$}}
%{\bf Florent NGUEMA NDONG$^1$}}
\\[2mm]
 $^{1}$ Universit\'e des Sciences et Techniques de Masuku 
 \\
Email: \href{mailto:florentnn@yahoo.fr}{\color{blue}{\underline{\smash{florentnn@yahoo.fr}}}}\\[1mm]
$^2$ Université de Poitiers\\
Email: \href{anne.bertrand@math.univ-poitiers.fr}{\color{blue}{\underline{\smash{anne.bertrand@math.univ-poitiers.fr}}}}\\[1mm]
%$^{3}$ Universit\'e des Sciences et Techniques de Masuku 
 \\
%Email: \href{mailto:florentnn@yahoo.fr}{\color{blue}{\underline{\smash{florentnn@yahoo.fr}}}}\\[1mm]
\lin{5}{7}

 {  %\fontfamily{Cambria}\selectfont
 {\noindent{\large\bf\color{myaqua} Abstract}{\bf \\[3mm]
 \textup{
  Let $ \beta $ be a real number smaller than -1. Endowed with the shift, the negative $\beta$-shift is a symbolic dynamical system that is coded only under certain conditions, but in all cases, we show that the support of the measure of maximal entropy is a free monoid generated by a positive recurrent code. One of the differences between the positive and negative $\beta$-shift is the existence of gaps in the system for certain negative values of $ \beta $. These are intervals of negative $\beta$-representations (cylinders) that are null sets with respect to the measure of maximal entropy, which is a Champernowne measure.
 }}}
 \\[4mm]
 {\noindent{\large\bf\color{myaqua} Keywords}{\bf \\[3mm]
  negative bases; $\beta$-expansions; coded systems; transitivity; intrinsic ergodicity 
}}
\\[4mm]
 {\noindent{\large\bf\color{myaqua} 2010 MSC}{\bf \\[3mm] 11K16; 11B05; 37A05; 37A25; 37B10
}}}
\lin{3}{1}

\renewcommand{\headrulewidth}{0.5pt}
\renewcommand{\footrulewidth}{0pt}

\newcounter{itemref}

 \pagestyle{fancy}
 \fancyfoot{}
 \fancyhead{} % clear all header and footer fields
 \fancyhf{}
% \fancyhead[RO]{\leavevmode \put(-105,0){\color{myaqua}F. NGUEMA NDONG} \boxx{15}{-10}{10}{50}{15} }
 \fancyfoot[C]{\leavevmode
 %\put(0,0){\color{lightaqua}\circle*{34}}
 %\put(0,0){\color{myaqua}\circle{34}}
 \put(-2.5,-3){\color{myaqua}\thepage}}

 \renewcommand{\headrule}{\hbox to\headwidth{\color{myaqua}\leaders\hrule height \headrulewidth\hfill}}

 \let\thefootnote\relax\footnotetext{*This work was partially supported by the project Mathematics in Gabon for Africa (ERC).}
 \let\thefootnote\relax\footnotetext{*This work was partially supported by the project International Emerging Actions (CNRS)}
%\addtocounter{footnote}{0} 
%\footnotetext{The research of authors on this project has received funding from the Centre National de Recherche Scientifique (CNRS).}

 \section{Introduction}

%{ %\fontfamily{times}
\selectfont
 \noindent 
 
Let $ \beta $ be a real number such that $ |\beta|>1$. Since the seminal paper \cite{MR0097374} on expansions of numbers in non-integer positive bases, many mathematicians have become interested in the properties of the $\beta$-shift for $ \beta > 1 $. For instance, in \cite{MR939059} and \cite{hofbauer1978beta}, the authors established the intrinsic ergodicity. The condition ensuring the specification property is given in \cite{MR545668}. The concept of intrinsic ergodicity was first explored by W. Parry in \cite{MR0142719}: it concerns the uniqueness of the measure of maximal entropy, and the study of the impossibility of decomposing the system into several invariant non-negligible subsets with respect to this measure. Recall that the set of invariant measures on a dynamical system is nonempty.

In \cite{Vitt}, Vittorio Grünwald introduced the negabinary system. It is a non-standard positional numeral system with the unusual property both that negative and positive numbers can be represented without a sign bit. It was used in the experimental Polish computers SKRZAT 1 (see for example \cite{FIETT19601250}) and BINEG in 1950. In 2009, S. Ito and T. Sadahiro, in \cite{MR2534912} extended the notion of signless numbers representation to all negative bases. They determined the unique measure of maximal entropy of $ I_{\beta}=[\frac{\beta}{1-\beta}, \frac{1}{1-\beta})$ endowed with the negative $\beta$-transformation.
The theory of representation of numbers highlights the connection between the symbolic dynamics and analytic number theory. 
Over the past two decades, various papers have investigated to the similarities and differences between the positive and negative $ \beta$-transformations, and several fundamental properties of strings associated with negative $\beta$-expansions have been established. See, for instance \cite{DBLP:journals/ijfcs/CharlierS18}, \cite{NguemaNdong20161}, \cite{NguemaNdong2019}, \cite{MR2974214}. 
This paper fits within this line of research.
Endowed with the shift, the positive and negative $\beta$-shifts share many properties. On the other hand, they also differ in several important aspects. For instance, the positive $\beta$-shift is coded for every value of $ \beta$ greater than 1. In contrast, the negative $\beta$-shift is coded if and only if $ \beta$ is smaller than or equal to $ -\frac{1+\sqrt{5}}{2} $ and the $ \beta$-expansion of the left end-point of the domain of the negative $\beta$-transformation is not periodic with odd period. 
 Moreover, a striking phenomenon distinguishes positive from negative $\beta$-shifts: the existence of gaps in the negative case. These gaps correspond to subsets (cylinders) of the negative $\beta$-shift that are negligible with respect to the maximal entropy measure. Such cylinders are carried by intransitive words (see  Theorems \ref{InW1} and \ref{InW2}). In the domain of the negative $\beta$-transformation (which we will denote by $I_{\beta}$), this phenomenon was closely investigated in \cite{MR2974214}. Moreover, the authors give a very simple and interesting proof of the uniqueness of the measure of maximal entropy on the negative $\beta$-transformation (see \cite{MR2974214} Corollary 2.4). This result suggests that the one-sided negative $\beta$-shift is intrinsically ergodic. Also, when $ \beta $ is smaller than $-\frac{1+\sqrt{5}}{2} $,  M. Shinoda and K. Yamamoto proved in \cite{article} that factors of the negative $\beta$-shift are intrinsically ergodic.

In this paper, we describe the unique measure of maximal entropy on the negative $\beta$-shift. We show that it is the Champernowne measure associated with a positive recurrent prefix code  (see Theorem \ref{th3}) and we give an explicit description of this code. This description allows us to understand the factors of the negative $\beta$-shift for $ -\frac{1+\sqrt{5}}{2} \leq \beta < -1 $.
Indeed, in this case, the negative $ \beta$-shift (denoted by $S_{\beta}$) is not transitive. But, the measure of maximal entropy is supported by a transitive subset strictly included in the system and obtained as image of a coded $x$-shift under a morphism, with $ x \leq - \frac{1+\sqrt{5}}{2}$ (see Theorem \ref{Th tau}).
We describe this support by identifying some transitive subshifts contained in $ S_{\beta} $ with  associated codes (prefix or suffix). To construct the unique measure of maximal entropy,
our approach mainly relies on well-known results on coded systems due to G. Hansel and F. Blanchard in \cite{MR858689}, as well as A. Bertrand in \cite{MR939059}. 

We highlight a surprising phenomenon which establishes a link between $\beta$-representations with $ \beta $ smaller than $ - \frac{1+\sqrt{5}}{2} $ and those for which $ \beta $ is taken in $ ] -\frac{1+\sqrt{5}}{2}, -1 [ $. There is a decreasing sequence of real numbers $ (\gamma_n)_{n \geq 0 } $ with limit 1 and $ \gamma_0 = \frac{1+\sqrt{5}}{2} $ and we define the map from $ \N $ to $\{0, 1\}^* $ such that $ \phi(k) = 1(00)^k$. We show that there is a bijection between $ ] -\infty, -\gamma_0 ]$ and each interval $ ]-\gamma_n, -\gamma_{n+1} ] $. For any $ \beta \in ]-\infty, -\gamma_0 ] $ there is a unique $ \beta_n \in ] - \gamma_n, - \gamma_{n+1} ]$ such that $ d(l_{\beta_n}, \beta_n) = \phi^{n+1}(d(l_{\beta}, \beta)) $. The $ \beta$-shift $ S_{\beta}$ is coded by a language $\C_{\beta}$ and the measure of maximal entropy on the two-side $\beta$-shift is carried by $ \C_{\beta}^{\Z}$. The measure of maximal entropy on the two-sided $ \beta_n$-shift is $ \phi^{n+1}(\C_{\beta})^{\Z}$. That of the right one-sided $ \beta_n$-shift is carried by $ \bigcup\limits_{k \geq 0} \sigma^k(\phi^{n+1}(S_{\beta})) $. As a consequence, if $\C_{\beta}$ codes $ S_{\beta}$, then $\phi^{n+1}(\C_{\beta})$ codes the support of the measure of maximal entropy on $S_{\beta_n}$.

%Note that, if it can be deduced from the intrinsic ergodicity of the negative $\beta$-transformation that of the negative $\beta$-shift, the existence of the positive recurrent code in the system offers another way to prove the uniqueness the measure of maximal entropy.

  The contents of this paper are as follows. We start our study by generalities on symbolic dynamical system. More precisely, we begin with a brief overview of coded systems, the notion of tower of a prefix code (introduced by G. Hansel and F. Blanchard in \cite{MR858689}) and $\beta$-shift. 
 The second part of the paper is devoted to the explicit description of codes and the ergodic measures attached to the negative $\beta$-shifts, and the surprising links between $ ]-\infty, -\gamma_0]$ and each interval $]-\gamma_{n}, -\gamma_{n+1}]$. We start this section by recalling the codes of the possible supports of the measure of maximal entropy. Next, we characterize intransitive words.

\section{Generality}

In this section, we briefly recall several facts about coded systems and the representation of numbers in non-integer real bases. 

\subsection{Coded System}

Let $\A$ be an alphabet, $(S,T)$ a symbolic dynamical system on $ \A$ and $ L_S$ the associated language. In the following, we denote by $ \A^*$ the free monoid generated by $ \A $, and $ \A^{+} = \A^{*}\backslash\{ \varepsilon \} $, where $ \varepsilon $ is the empty word, $ \R$, $\Z$ and $\N $ denote respectively, the set of real numbers, the set of integers and the set of nonnegative integers; $ \N^{\times} = \N\backslash\{0\} $, $\R^{\times} = \R\backslash\{0\}$ and $ \Z^{\times} = \Z\backslash\{0\}$. We recall some definitions given in \cite{MR858689}.

\begin{defi}
A language $L$ is said to be transitive if for any pair of words $(u,v)$ of $L^2=L\times L$, there exists $w$ in $ \mathcal{A}^{*} $ such that  $uwv$ belongs to $L$.
\end{defi}

A symbolic dynamical system is said to be transitive if its language is transitive. More generally, given a topological dynamical system $ (S,T)$, for all open sets $U$ and $ V$ of $ S $, we have $ U \cap T^{-n}V \neq \varnothing $ for some $ n$ in $\Z$. This is equivalent to saying that the orbit $ \bigcup \limits_{n \in \Z} T^n U $ of all non-empty open set $ U $ of $S$ is dense in $ S $ (see \cite{MR0352411}).

\begin{defi}
A code $Y$ on $\A$ is a language such that, for any equality 
\begin{equation}
x_1x_2\cdots x_n = y_1y_2 \cdots y_k,
\end{equation}
for any $ x_i$, $y_j \in Y$, we have $ x_i=y_i $ and $ k = n$. 
\end{defi}

A prefix (resp. suffix) code is a language $ \C $ of $\A^{+}$ for which no word is a prefix (resp. suffix) of another. That is,
\begin{equation}
 \forall u, v \in \C, u = vw \Rightarrow u = v \text{ and } w = \varepsilon.
\end{equation}

\begin{defi}
A language $ L $ is said to be coded by $ X $ if $ L$ is the set of factors of $ X^*$, the free monoid generated by $ X $. A system $ S $ is coded if its language is coded.
\end{defi}

One endows the alphabet $\A$ with the dyadic metric $ d $ defined by: for both sequences $ x= (x_i)_{i \in \Z} $ and $ y= (y_i)_{i \in \Z} $ of $ \A^{\Z}$,
\begin{equation}
 d(x, y) = \sum \limits_{n \in \Z} 2^{-\vert n \vert} d(x_n , y_n) \label{diadic 1}
\end{equation}
 with
\begin{equation}
 d(x_n, y_n)= \begin{cases}
               0 &\text{ if $ x_n = y_n $ } \\
               1 &\text{ if $ x_n \neq y_n $ }.
              \end{cases}\label{diadic 2}
\end{equation}
Let $ S $, endowed with the metric $ d $, be a symbolic dynamical system over $ \A$ coded by $ X$. If $ S $ is a two-sided system, then, it is the closure of $ X^{\Z}$. That is, the strings which belong to $ S $ are infinite concatenations of words of $ X$ or limits of sequences of strings of $ X^{\Z}$. However, if $ S $ is a right one-sided symbolic dynamical system, then $ x \in S$ if
$ x = u x_1x_2\cdots $, where $ (x_i)_{i \geq 1}$ belongs to the closure of $ X^{\Z}$ and $ u $ is the empty word or a suffix of a word of $ X$.

Let $ L $ be a language on an alphabet $ \A $. The radius $ \rho_L $ of the power series $ \sum\limits_{n \geq 1} card (L\cap \A^n)z^n $ is called the radius of convergence of $ L $. Since $ L^* $ is the free monoid generated by $ L $, $ \rho_{L^*} $ denotes the radius of convergence of the power series $ \sum\limits_{n \geq 1} card (L^*\cap \A^n)z^n $. 

\begin{defi}
Let $ S $ be a symbolic dynamical system and $ a_1 a_2 \cdots a_k \in L_S$. We denote by $ _{m}[a_1\cdots a_k ]$, the set of sequences $ (x_i)_{i \in \Z}$ of $ S $ such that 
\begin{equation}
x_m x_{m+1} \cdots x_{m+k-1} = a_1 a_2 \cdots a_k.
\end{equation}
This set is called cylinder carried by $ a_1a_2\cdots a_k $ at $ m $. For $ m = 0 $, in order to simplify the notations, we will denote  $_{0}[x] $ by $ [x] $.
\end{defi}

Consider a symbolic dynamical system $ S $ on an alphabet $ \A $ and $ x \in L_S $. The length of a word $ x $, denoted by $ l(x) $, is the number of letters of $ x $. 
\begin{equation}
x = x_1 x_2 \cdots x_k, \text{ with $ x_i \in \A $ }\Rightarrow l(x) = k.
\end{equation}

\begin{defi}
A prefix code $ \C $ is said to be positive recurrent if 
\begin{equation}
\sum\limits_{x \in \C} \rho_{\C^{*}}^{ l( x ) } = 1 \text{ and } \sum\limits_{ x \in \C} l( x)  \rho_{\C^{*}}^{l( x )}<+\infty.
\end{equation}

If we set $ \rho_{\C^{*}} = \frac{1}{\beta}$, then $ 1 = \sum\limits_{n \geq 1} \frac{c_n}{\beta^n} $ and $ \sum\limits_{n \geq 1}\frac{nc_n}{\beta^n} < +\infty $ where $ c_n $ counts the number of words of length $ n $ in $ \C $ and $ \C^{*} $ the free monoid generated by $ \C $.
\end{defi} 

\begin{rem}\label{Nrem1}
From Proposition 2.15 of \cite{MR858689}, if $\mu$ is an invariant measure on $\C^{\Z}$ endowed with the shift and $h(\mu)$ is its entropy, then
 \begin{itemize}
  \item[(1)] we have
  \begin{equation}
   h(\mu) \leq - l(\C, \mu)\log \rho_{\C^{*}}
   \label{(1)}
  \end{equation}
  where
  \begin{equation}
  l(\C, \mu) = \sum\limits_{ x \in \C }l(x)  \mu([x]).
  \end{equation}
\item[(2)] In \eqref{(1)}, the equality holds if both following conditions are satisfied:
  \item[(a)] $ \sum \limits_{x \in \C} \rho^{l( x )}_{\C^{*}} = 1 $,
  \item[(b)] $ \mu $ is a probability of Bernoulli on $\C^{*}$ defined by
  $ \mu([x])=\rho^{l(x)}_{\C^{*}}$, $ x \in \C $.
 \end{itemize}
\end{rem}
 
\begin{defi}
A measurable topological dynamical system $ (S, m, g) $ is said to be \textit{ergodic} if for any measurable $g$-invariant set $ B \subset S $, we have $m(B) = 0$ or $ m(B) = 1$. One also says that $ m $ is an ergodic measure with respect to $ g $ or that $ g $ is ergodic with respect to $ m $. It is said to be mixing if for all measurable sets $ A $ and $ B $,
 \begin{equation}
  \lim \limits_{n \rightarrow +\infty} m(g^{-n}(A) \cap B) = m(A)m(B).
 \end{equation}
\end{defi}
Let $ \A$ be an alphabet. In what follows, $ \sigma$ denotes the shift: the map from $\A^*$ to $ \A^*$ such that $ \sigma((x_n)_{n\in\Z})=(x_{n+1})_{n\in \Z}$.

\subsection{Tower associated to a prefix code}

More details on the notion of the tower of a prefix code can be found in \cite{MR858689}.
Let $ \Omega $ be the subset of $ \C^{\Z}\times \N $ such that:
\begin{equation}
 ((x_n)_{n \in \Z}, i) \in \Omega \Rightarrow 1 \leq i \leq l( x_0 ).
\end{equation}
We can identify $ (x_n)_{n \in \Z} $ with an element $ x $ of $ \C^{\Z}$ which is a concatenation of words $ x_i $ of the code $\C$.

We define a map $ T $ from $ \Omega $ into itself by:
\begin{equation}
 T((x_n)_{n \in \Z}, i) = \begin{cases}
                           ((x_n)_{n \in \Z}, i+1)   & \text{ if $ i < l( x_0 ) $ } \\
                           ((x_{n+1})_{n \in \Z}, 1) & \text{ if $ i = l( x_0 ). $ }
                          \end{cases}
\end{equation}
The pair $ (\Omega, T) $ is called \textit{the tower associated} to $ \C $. We have the two important following facts (see \cite{MR858689}).
 \begin{itemize}
  \item[•] When $ \overline{\mu} $ ranges through $\M_T(\Omega)$, $\sup_{\overline{\mu}}h(\overline{\mu}) = -\log \rho_{\C^{*}}$. 
  
  \item[•] There exists one and only one invariant probability $ \mu $ on $ \Omega $ such that $ h(\overline{\mu})=-\log \rho_{\C^{*}} $ if and only if $\C$ is positive recurrent. In this case, $ \overline{\mu} $ is the unique invariant probability (and thus ergodic) on $ \Omega $ inducing on $ \C^{\Z} $ a probability $\mu$ of Bernoulli defined by:
   \begin{equation}
   \mu ( [ x ] )= \rho^{l(x)}_{\C^{*}}, \text{ $ x \in \C $ }.
  \end{equation}
 \end{itemize}
The first statement is usually referred to as the variational principle. When the dynamical system is coded by a recurrent prefix code $\C$, the set of periodic points is dense. We say that two words $ u $ and $ v $ of the language $L$ are in the same class of the syntactic monoid if for all pair of words $(a,b)$,
\begin{equation}
 aub \in L \Longleftrightarrow a v b \in L. 
\end{equation}
A symbolic dynamical system $S$ is said to be rational (or sofic) if the number of classes of language associated to is finite. In this case, it is coded by a positive recurrent prefix code.
The measure $ \overset{\bullet}{\mu} $ on the coded system induced by the measure on the tower is particularly simple: if $ u_1u_2 \cdots u_r $ is a word of $ L $,
\begin{equation}
 \overset{\bullet}{\mu}([u_1u_2 \cdots u_r]) = \sum\limits_{ m } \mu([m]) = \dfrac{1}{\sum\limits_{n \geq 1}\frac{n.c_n}{\beta^n}}
 \sum \limits_{ m } \frac{1}{\beta^{l( m )}}
\end{equation}
where in the sums, $ m $ is all words of the form $ m = a u_1u_2 \cdots u_r b $; $ a $ is a proper prefix of a word of the code and $ b $ is a proper suffix, or the empty word.

Each class of the syntactic monoid is associated to a positive constant $ \lambda $ such that if $ u $ belongs to this class, the cylinder $[u]$ has a measure
\begin{equation}
 \overset{\bullet}{\mu}([u]) = \frac{\lambda}{\beta^{l(u )}}.
\end{equation}
So, the measure of a cylinder depends on the class of its support (the word $ u $) and its length. If there exist two positive numbers $ \epsilon $ and $ M $, with  
 $ 0 < \epsilon < \lambda < M $ for all classes, the measure is said to be homogeneous. The ratio between measures of cylinders of the same length is controlled by $ \frac{\epsilon}{M} $ and $ \frac{M}{\epsilon} $. If the greatest common divisor ($\gcd$) of lengths of words of the code is 1, then the measure is mixing. In this case, it is said to be a Champernowne measure.

\subsection{Beta-shift }

In this section, we provide definitions that unify the cases of $ \beta< -1$ and $\beta>1$. 
Consider a real number $ \beta $ with modulus greater than 1. Let us approach the question of the representation of numbers using powers of $ \beta $ from a more general point of view. We set 
\begin{equation}
l_{\beta} = \begin{cases}
          0 &\text{ if $ \beta > 1 $ } \\
          \frac{\beta}{1-\beta} &\text{ if $ \beta < -1 $ }
           \end{cases} \text{ and $ r_{\beta} = l_{\beta} + 1 $. }
\end{equation} 
\subsubsection{Beta-transformation and beta-expansions}
Define the map $ T_{\beta}$ from $I_{\beta} = [l_{\beta}, r_{\beta} )$ into  itself by  
\begin{equation}
  T_{\beta}(x) = \beta x - \lfloor \beta x -l_{\beta} \rfloor. \label{T}
\end{equation}

The expansion in base $ \beta $ of a real $ x $ (denoted by $ d(x, \beta )$) is given by the following algorithm.
We find the smallest nonnegative integer $ n $ for which $ \frac{x}{\beta^n} \in I_{\beta} $. The $\beta$-expansion is given by the sequence $ d(x, \beta) =x_{-n+1}\cdots x_0 \bullet x_1x_2 \cdots $ such that
\begin{equation}
  x_{-n+i} = \lfloor \beta T_{\beta}^{i-1}\left(\frac{x}{\beta^n}\right)-l_{\beta} \rfloor,\text{ $ i \geq 1 $ }.
\end{equation}

\subsubsection{Characteristic sequences}
Now we define the sequence controlling the $\beta$-shift.
%In what follows, $ (d_i)_{i \geq 1}$ denotes the $\beta$-expansion of $ l_{\beta} $.
Let $ (r_i^*)_{i \geq 1 } $ be the sequence of digits such that:
\begin{equation}
r^*_1 = \lfloor \beta r_{\beta} - l_{\beta} \rfloor
\end{equation}
and for all integers $ i \geq 2 $,
\begin{equation}
\begin{aligned}
r^*_i &= \lfloor \beta T_{\beta}^{i-2}(\beta r_{\beta}-r^*_1)-l_{\beta} \rfloor\\
    &=\lfloor \beta^i r_{\beta}- \sum\limits_{k=1}^{i-1} r^*_k \beta^{i-k}-l_{\beta} \rfloor.
\end{aligned}
\end{equation}
Consider an alphabet $ \A $ endowed with an order $ \prec_{\delta}$ such that for all sequences of digits $(x_i)_{i \geq 1} $ and $(y_i)_{i \geq 1}$ over $ \A $, 
\begin{equation}
(x_i)_{i \geq 1} \prec_{\delta} (y_i)_{i \geq 1} \Leftrightarrow \exists k \in \N^{\times}| \hspace{0.2 cm} x_i=y_i \hspace{0.2cm}\forall i < k, \hspace{0.5cm}\delta^k(x_k-y_k)< 0
\end{equation}
where $ \delta $ is the sign of $ \beta $. Furthermore,
\begin{equation*}
(x_i)_{i \geq 1} \preceq_{\delta} (y_i)_{i \geq 1} \Leftrightarrow (x_i)_{i \geq 1} \prec_{\delta} (y_i)_{i \geq 1} \text{ or } (x_i)_{i\geq 1} = (y_i)_{i \geq 1}.
\end{equation*}
\begin{itemize}
\item If $ \beta < -1 $, $ \prec_{\delta} $ is the \textit{alternating order} (see for example \cite{MR2534912}, \cite{NguemaNdong20161}, \cite{NguemaNdong2019}).
\item If $ \beta> 1$, $ \prec_{\delta} $ is the \textit{classical lexicographic order} on words. 
\end{itemize}

Let $ d(l_\beta, \beta) = \bullet d^*_1 d^*_2 \cdots $. If $ \beta > 1 $, then for any integer $ i $, $ d_i^* = 0 $. 
Let $ (r_i)_{i \geq 1}$ be the sequence of digits such that
\begin{equation}
 (r_{i})_{i \geq 1} = \begin{cases}
                         \overline{(r_1^*, \cdots, r_{n-1}^*, r_{n}^*-1)} &\text{ if $ (r_i^*)_{i \geq 1} = (r_1^*, \cdots, r_{n}^*,d^*_1,d^*_2,\cdots ) $, $\beta>1$ }\\
                         \overline{(r_1^*, \cdots, r_{n-1}^*, r_{n}^*-1)} &\text{ if $ (r_i^*)_{i \geq 1} = (r_1^*, \cdots, r_{n}^*,d^*_1,d^*_2,\cdots ) $, $\beta<-1$ and $ n $ even }\\
                                  (r_i^*)_{i \geq 1} &\text{ otherwise } 
                      \end{cases}
\end{equation} 
and
\begin{equation}
 (d_{i})_{i \geq 1} = \begin{cases}
                         \overline{(d^*_1, \cdots, d^*_{2n-2}, d^*_{2n-1}-1, 0)} &\text{ if $\beta< -1$, $ (d^*_i)_{i \geq 1} = \overline{(d^*_1, \cdots, d^*_{2n-1} )} $ }\\
                         (d^*_i)_{i \geq 1} &\text{ otherwise }
                         \end{cases}\label{(D)}
\end{equation}      
where $ \overline{t} $ stands for infinite repetition of the word $ t $.
Indeed, 
\begin{equation}
(d_i)_{i \geq 1} = \lim\limits_{x \rightarrow l_{\beta}^+}d(x, \beta) \text{ and } (r_i)_{i \geq 1} = \lim\limits_{x \rightarrow r_{\beta}^-}d(x, \beta) \label{lbrb}
\end{equation}
and we set:

\begin{equation}
d^{*}(l_{\beta}, \beta) = (d_i)_{i \geq 1}  \text{ and } d^{*}(r_{\beta}, \beta)=(r_i)_{i \geq 1} . \label{lbrbcarac}
\end{equation}

 \textit{Expansions in base} $ \beta $ (or $ \beta$-expansions, $ |\beta|> 1 $) of real numbers are governed by the sequences $(d_i)_{i \geq 1}$ and $(r_i)_{i \geq 1}$.
A sequence $(x_i)_{i \geq 1} $ is the $ \beta $-expansion of a real $ x $ if and only if 
\begin{equation}
(d^*_i)_{i \geq 1} \preceq_{\delta}(x_{i+n})_{i \geq 1} \prec_{\delta}(r_i)_{i \geq 1}, \text{ $ \forall n \in \N $}.
\end{equation}
Particularly
\begin{equation}
 (d^*_i)_{i \geq 1} \preceq_{\delta}(d^*_{i+n})_{i \geq 1},  (r^*_{i+n})_{i \geq 1} \prec_{\delta}(r_i)_{i \geq 1}, \text{ $ \forall n \in \N $}. \label{lyn 1}
\end{equation}

\begin{defi}

In what follows, for $ \beta < -1$, $ d^{*}(l_{\beta}, \beta)$ is called \textit{characteristic sequence} associated with $ \beta$. 

If $\beta> 1$, the characteristic sequence refers to $ d^{*}(r_{\beta}, \beta) = d^{*}(1, \beta)$. 
\end{defi}

One determines $ (d_i^*)_{i \geq 1}$ with the transformation $T_{\beta}$ or by using the map $ f_{\beta} $  on words defined by:
\begin{equation}
f_{\beta} ((x_i)_{i \in \Z}) = \sum\limits_{ k \in \Z} x_k \beta^{-k}.
\end{equation}
If $d(x, \beta ) = x_{-n}x_{-n+1}\cdots x_{-1}x_0\bullet x_1x_2 \cdots$, then
\begin{equation}
x = f_{\beta} ((x_i)_{i \geq -n}). \label{lyn 2}
\end{equation}
Thus 
\begin{equation*}
r_{\beta} = f_{\beta} (d^{*}(r_{\beta}, \beta)) = f_{\beta}(d(r_{\beta},\beta)),
\end{equation*}
and
\begin{equation*}
l_{\beta} = f_{\beta} (d^{*}(l_{\beta}, \beta)) = f_{\beta}(d(l_{\beta},\beta)). 
 \end{equation*} 
 It is known that $f_{\beta}$ is an increasing map in the meaning of the order $ \preceq_{\delta}$. As consequence, 
 \begin{equation*}
 d^{*}(l_{\beta}, \beta) = \sup \{ (x_i)_{i \geq 1} : f_{\beta}((x_i)_{i \geq 1}) = l_{\beta} \text{ and } (x_i)_{i \geq 1 } \preceq_{\delta} (x_{i+n})_{i \geq 1} \text{ $ \forall n \in \N $} \}
\end{equation*} 
and 
 \begin{equation*}
 d^{*}(r_{\beta}, \beta) = \inf \{ (x_i)_{i \geq 1} : f_{\beta}((x_i)_{i \geq 1}) = r_{\beta} \text{ and } (x_{i+n})_{i \geq 1 } \preceq_{\delta} (x_{i})_{i \geq 1} \text{ $ \forall n \in \N $} \}.
\end{equation*} 
 
 Note that, the fact that a sequence $ (x_n)_{n \geq 1} $ satisfies conditions \eqref{lyn 1} and \eqref{lyn 2} does not mean that $ (x_n)_{n \geq 1}$ is the $ x $-expansion of $l_x$, for some $ x $. In negative bases cases, the conditions for a sequence to be the $ x $-expansion of $l_x$, for some $ x $ are stated in \cite{NguemaNdong20161}.

\subsubsection{Definition of the beta-shift}

\begin{defi}
The $\beta$-shift $\tilde{S}_{\beta}$  is the closure of the set of $\beta$-expansions.
\begin{equation}
\tilde{S}_{\beta} = \{ x_k x_{k+1} \cdots x_0 \bullet x_1\cdots \vert (d_i^*)_{i\geq 1} \preceq_{\delta}(x_i)_{i \geq m} \preceq_{\delta}(r_{i})_{i \geq 1}, \forall m\geq k, \forall k \}\label{sb}.
\end{equation}
We also define the corrected $\beta$-shift $ S_{\beta} $ as follows:
\begin{equation}
S_{\beta} = \{ x_k x_{k+1} \cdots x_0 \bullet x_1\cdots \vert (d_i)_{i\geq 1} \preceq_{\delta}(x_i)_{i \geq m} \preceq_{\delta}(r_{i})_{i \geq 1}, \forall m\geq k, \forall k \}.
\end{equation}
\end{defi}
Observe that $ S_{\beta}$ and $ \tilde{S}_{\beta}$, endowed with the shift $\sigma$ have the same entropy. Moreover, each real number has a representation in $ S_{\beta}$. Both bounds $(r_i)_{i \geq 1}$ and $ (d_i)_{i \geq 1}$ decide whether a digit string belongs to $ S_{\beta}$ or not. Thus it is more natural and convenient to use $S_{\beta}$ as $\beta$-shift instead of $\tilde{S}_{\beta}$. Indeed, $ S_{\beta}$ is the smallest $ \sigma$-invariant subshift with entropy $ \log|\beta|$ included in $ \tilde{S}_{\beta}$ and in which each real number has at least one $\beta$-representation. In the rest of the paper, $ L_{\beta} $ denotes the language of $ S_{\beta} $.

\begin{defi}
In what follows, sequences of $ S_{\beta} $ (or words of its language) are said to be admissible.
\end{defi}
When $ \beta $ is strictly smaller than -1, thanks to the map $ f_{\beta}$, it is proved that there are real numbers with several $ \beta$-representations, but the $ \beta$-expansion is unique. Indeed, if we set:
\begin{equation*}
\begin{cases}
x&= x_1 \cdots x_{k-1}(x_k+1)d_1d_2 \cdots , \\
y&= x_1 \cdots x_{k-1}x_k0d_1d_2 \cdots ,
\end{cases}
\end{equation*}
we have $ f_{\beta}(x) = f_{\beta}(y)$.
But $ x_1 \cdots x_{k-1}(x_k+1)d_1d_2 \cdots $ is the $ \beta$-expansion of $f_{\beta}(x) = f_{\beta}(y) $.

\begin{rem}\label{rmbta}
Let $ X $ be a symbolic dynamical system, $ \log t $ its entropy. Then, $ \frac{1}{t} $ is the smallest pole in modulus of the formal power series $ \sum\limits_{ x\in L_X} z^{l(x)} $.
\end{rem}

Denote by $ H_n $, the number of admissible words of length $ n $ in $ S_{\beta}$. From \cite{NguemaNdong20161},
\begin{equation}
 H_n = \sum\limits_{k = 1}^n\delta^n(d_{k-1}-d_k)H_{n-k}+1.
\end{equation}
Moreover, $ \lim\limits_{n \rightarrow + \infty} \frac{H_{n+1}}{H_n} = |\beta|$ and $ \log |\beta|$ is the topological entropy of the system. 

\begin{defi}\label{InS}
Let $ \beta$ be a real number smaller than $ -1 $ and $ (d_i)_{i \geq 1} $ its characteristic sequence. In what follows, the block $ d_1 \cdots d_k $ will be called \textit{initial segment of length $ k $}. Moreover, in the sequel, the $\beta$-expansion of $ l_{\beta}$ will denote the sequence $(d_i)_{i \geq 1}$ defined in \ref{(D)} (the corrected $\beta$-expansion).
\end{defi}

In what follows, $ S_{\beta}^r $ denotes the right one-sided negative $\beta$-shift, that is, the closure of $\beta$-expansions of real numbers of $I_{\beta} = [l_{\beta}, r_{\beta} )$. 

\section{A surprising phenomenon}

A surprising phenomenon occurs and establishes a link between $ ]-\infty, -\gamma_0 ]$ and the sequence of intervals $ ]-\gamma_n, -\gamma_{n+1}] $. %$\beta$-expansions of various intervals of $\beta$.

\subsection{The map $ \phi $}

Let $\phi$ be the morphism from $\N$ to $\{0, 1\}^*$ defined by:
\begin{equation}
\begin{aligned}
 \phi :  &\N && \rightarrow \{0, 1\}^* \\
         & k &&\mapsto \phi(k)=1(00)^k.
\end{aligned}
\end{equation}

We set $u_{-1}= 0$ and for all $ n$ in $ \N$, $ u_n = \phi^n(1) $, $ v_n = u_{n-1}u_{n-1} = \phi^n(00)$ and $ \gamma_n $ denotes the real number such that $ \frac{1}{\gamma_n} $ is the smallest real satisfying:

\begin{equation}
1=\frac{1}{\gamma_n^{l(u_n)}}+\frac{1}{\gamma_n^{l(v_n)}}.
\end{equation}
Indeed, $\gamma_n$ is the algebraic integer of the polynomial $ X^{l_n}-X-1 $, where $ l_n = \max(l(u_n),l(v_n))$ (see the proof of Proposition 5 in \cite{NguemaNdong20161}). Note that $ \gamma_0 $ is the golden ratio $ \frac{1+\sqrt{5}}{2} $. Moreover,
\begin{equation}
d(l_{-\gamma_0}, -\gamma_0) = 1 (00)^{\infty}= \lim\limits_{k\rightarrow +\infty} \phi(k).
\end{equation}
More generally (see \cite{MR2974214}), 
\begin{equation}
d(l_{-\gamma_n}, -\gamma_n) = u_n (u_{n-1}u_{n-1})^{\infty}.
\end{equation}

\begin{exple}
Let $ \beta $ be the real number smaller than -2 such that $d(l_{\beta}, \beta) = 2\overline{1}=2(1)^{\infty}$. Then:

\begin{equation*}
\phi (2\overline{1}) = 10000(100)^{\infty}.
\end{equation*}
\end{exple}

\begin{prop}\label{Tau Inc}

 In the meaning of the alternating order, $\phi $ is an increasing map.
\end{prop}

\begin{preu}\textbf{of Proposition \ref{Tau Inc}}

Let $(x_i)_{i \geq 1}$ and $(y_i)_{i \geq 1}$ be two sequences of nonnegative integers such that $(x_i)_{i \geq 1}\prec (y_i)_{ \geq 1}$. There is $ k $ satisfying $x_1 \cdots x_{k-1} = y_1 \cdots y_{k-1}$ and $(-1)^k(x_k-y_k) < 0$. We have $ \phi ((x_i)_{i \geq 1}) = 1(00)^{x_1}1(00)^{x_2} \cdots$ and $\phi((y_i)_{i \geq 1}) = 1(00)^{y_1}1(00)^{y_2}\cdots$.

\begin{equation*}
 x_1 \cdots x_{k-1} = y_1 \cdots y_{k-1}\Rightarrow 1(00)^{x_1}\cdots 1(00)^{x_{k-1}} = 1(00)^{y_1}\cdots 1(00)^{y_{k-1}}.
\end{equation*}
And

\begin{equation*}
(-1)^k(x_k-y_k)<0 \Rightarrow \begin{cases}
                               1(00)^{x_k}1\prec 1(00)^{y_k}1&\text{if $k$ is odd } \\
                               1(00)^{y_k}1\prec 1(00)^{x_k}1&\text{if $k$ is even}.
                              \end{cases}
\end{equation*}
Using properties of the alternating order given in Section 1.1 of \cite{NguemaNdong2019}, we have

\begin{equation*}
 1(00)^{x_1}1\cdots1(00)^{x_{k-1}}1(00)^{x_k} \prec 1(00)^{y_1}1\cdots1(00)^{y_{k-1}}1(00)^{y_k},
\end{equation*}
and the result follows.
\end{preu}
\\
An immediate consequence of Proposition \ref{Tau Inc} is the following Corollary.

\begin{coro}\label{Cor Tau}
If $(d_i)_{i \geq 1}=d^{*}(l_{\beta}, \beta)$ is the characteristic sequence associated with  $\beta$ for some $\beta<-1$, then for any integer $n$, there exists $ x < -1$ such that $\phi^n((d_i)_{i \geq 1})=d^{*}(l_x, x)$ is the characteristic sequence associated with $x$.
\end{coro}

A sequence $(a_i)_{i \geq 1}$ is a $\beta$-expansion of $l_{\beta}$ for some $\beta<-1$ if it satisfies the following conditions (see \cite{Steiner2013}):

\begin{enumerate}[label=(\roman*), ref=\roman*]
\item\label{(a)} $ (a_i)_{i \geq 1} \preceq (a_{n+i})_{i \geq 1}$, for all integer $ n \in \N$;
\item\label{(b)} $ (a_i)_{i \geq 1} \prec \phi^{\infty}(1)$; 
\item\label{(c)} $ f_{\beta} ((a_i)_{i \geq 1}) = \frac{\beta}{1-\beta}$; 
\item\label{(d)} $ (a_i)_{i \geq 1} \not\in \{ a_1\cdots a_{k}, a_1\cdots a_{k-1}(a_k-1)0 \}^{\N} \setminus \{\overline{a_1\cdots a_k}\}$.

The condition \ref{(d)} means that if $(a_i)_{i\geq 1}$ belongs to  $ \{ a_1\cdots a_{k}, a_1\cdots a_{k-1}(a_k-1)0 \}^{\N}$, then it is equal to the periodic sequence $ \overline{a_1\cdots a_k}$. If $ k $ is odd, for the construction of the corrected $\beta$-shift $ S_{\beta}$, we choose as $\beta$-expansion of $l_{\beta}$, the sequence $\overline{a_1\cdots a_{k-1}(a_k-1)0}$. In the meaning of the alternating order, it is the largest sequence of the set $ \{ a_1\cdots a_{k}, a_1\cdots a_{k-1}(a_k-1)0 \}^{\N} $.

Let $\min(d_1\cdots d_k)$ and $ \max(d_1\cdots d_k)$ be respectively the minimum and the maximum of the interval of words $\{ d_1\cdots d_{k}, a_1\cdots d_{k-1}(d_k-1)0 \}^{\N} $, with $ d_k \neq 0$. We have:
\begin{equation}
\min(d_1\cdots d_k) = \begin{cases}
                       d_1\cdots d_{k-1}(d_k-1)0\overline{d_1\cdots d_k} &\text{ if $ k $ even} \\
                       d_1\cdots d_k \overline{d_1\cdots d_{k-1}(d_k-1)0} &\text{ if $k $ odd}
                       \end{cases}
\end{equation}
and
\begin{equation}
\max(d_1\cdots d_k) = \begin{cases}
                       \overline{d_1\cdots d_k} &\text{ if $ k $ even} \\
                       \overline{d_1\cdots d_{k-1}(d_k-1)0} &\text{ if $k $ odd}.
                       \end{cases}
\end{equation}
%If $(d_i)_{i \geq 1}$ satisfies \ref{Lyn 1} or \ref{Lyn 2}, then $ d(l_{\beta}, \beta) = \overline{d_1 \cdots d_k}$. Otherwise, $ d(l_{\beta}, \beta) = (d_i)_{i \geq 1}$. In fact, $ \overline{d_1\cdots d_k} $ is the expansion obtained by the greedy algorithm. But, $\max(d_1\cdots d_k )$ is the corrected $\beta$-expansion of $l_{\beta}$ used in the construction of the corrected negative $\beta$-shift $S_{\beta}$.
Thus, the characteristic sequence $d^{*}(l_{\beta}, \beta)=(d_i)_{i \geq 1}$ satisfies \ref{(a)}, \ref{(b)} and \ref{(c)}. 
%\begin{itemize}
\item\label{(e)} If $d^{*}(l_{\beta}, \beta) \in \{ a_1\cdots a_{k}, a_1\cdots a_{k-1}(a_k-1)0 \}^{\N}$, then $ d^{*}(l_{\beta}, \beta)=\max(a_1\cdots a_k)$.
\end{enumerate}

\begin{preu}{\textbf{Corollary \ref{Cor Tau}}}

Let $ \beta < -1 $ and $d^{*}(l_{\beta}, \beta)=(d_i)_{i \geq 1}$ be the characteristic sequence associated with $\beta$. Since $ \phi $ is increasing and $\phi^{\infty}(1)$ is invariant under $ \phi$, it follows that
\begin{equation*}
 (d_i)_{i \geq 1} \prec \phi^{\infty}(1) \Rightarrow \phi((d_i)_{i\geq 1})\prec \phi^{\infty} (1).
 \end{equation*}
 And thus \ref{(b)} is satisfied. 
 
 Moreover, for any integer $ n \geq 1$,  
\begin{equation*}
 (d_i)_{i \geq 1} \preceq (d_{i+n})_{i \geq 1} \Leftrightarrow 1(00)^{d_1}1(00)^{d_2}\cdots \preceq 1(00)^{d_n}1(00)^{d_{n+1}}\cdots.
 \end{equation*} 
 Indeed, $\phi((d_i)_{i \geq1})$ is smaller than all of its sub-words beginning with 1. Since the first letter of $\phi((d_i)_{i \geq1})$ is 1, it follows that it is less than all of its infinite sub-words beginning with 0 too. Then \ref{(a)} is satisfied. 
 
 It is proved in \cite{NguemaNdong20161} that the conditions \ref{(a)} and \ref{(b)} implies \ref{(c)}. Thus there is $ x< -1$ such that 
 \begin{equation*} f_{\beta}(\phi((d_i)_{i \geq1})) = \frac{x}{1-x}.
 \end{equation*} 
 %Thus, this sequence satisfies the following condition: 
%\begin{itemize}
%\item[•] $ (d_i)_{i \geq } \preceq (d_{i+n})_{i\geq 1} \prec (d_{i-1})_{i\geq 1} $, with $ d_0 = 0$;
%\item[•] 
%\end{itemize}
%Suppose $ \max(d_1\cdots d_k) \prec \phi^{\infty}(1)$, 
Note that $ \phi $ maps $\{ d_1\cdots d_{k}, d_1\cdots d_{k-1}(d_k-1)0 \}^{\N} $
 to a set of the same type. Indeed,
\begin{equation*}
\begin{aligned}
\phi(\min(d_1\cdots d_{2n}))&= 1(00)^{d_1}\cdots 1(00)^{d_{2n-1}}1(00)^{d_{2n}-1}1\overline{1(00)^{d_1}\cdots 1(00)^{d_{2n-1}}1(00)^{d_{2n}}}\\
                            &= \min (1(00)^{d_1}\cdots 1(00)^{d_{2n-1}}1(00)^{d_{2n}-1}1);
\end{aligned}
\end{equation*}
and
\begin{equation*}
\begin{aligned}
\phi(\min(d_1\cdots d_{2n+1}))&= 1(00)^{d_1}\cdots 1(00)^{d_{2n+1}}1\overline{1(00)^{d_1}\cdots 1(00)^{d_{2n}}1(00)^{d_{2n+1}-1}1}\\
                            &= \min (1(00)^{d_1}\cdots 1(00)^{d_{2n}}1(00)^{d_{2n+1}-1}1).
\end{aligned}
\end{equation*}
More precisely, 
\begin{equation*}
\phi (\min(d_1\cdots d_k)) = \min ( 1(00)^{d_1}\cdots 1(00)^{d_{k-1}}1(00)^{d_k-1}1 ).
\end{equation*}
 By analogy, we prove that 
 \begin{equation*}
 \phi (\max(d_1\cdots d_k)) = \max ( 1(00)^{d_1}\cdots 1(00)^{d_{k-1}}1(00)^{d_k-1}1 ). 
 \end{equation*}
 
Conversely, for all Lyndon word $(a_i)_{i \geq 1}$ such that there is $ p $ with $a_p\neq 0 $ and
\begin{equation*}
1(00)^{\infty} \prec \min (a_1\cdots a_p) \preceq (a_i)_{i \geq 1} \preceq \max(a_1\cdots a_k) \prec \phi^{\infty}(1),
\end{equation*}
We have $ a_1\cdots a_p = 1(00)^{k_1}\cdots 1(00)^{k_m}1 $, with $ p=2\sum\limits_{i=1}^{m}k_i+m+1$. So, 
\begin{equation*}
 \min(a_1\cdots a_p)=\phi(\min (k_1\cdots k_{m-1}(k_m+1)))
 \end{equation*}
  and 
  \begin{equation*} 
  \max(a_1\cdots a_p) = \max(\phi(k_1\cdots k_{m-1}(k_m+1))).
  \end{equation*} 
 Thus if $ \phi((d_i)_{i\geq 1}) \in \{a_1\cdots a_p, a_1\cdots a_{p-1}(a_p-1)0\}$, then $ (d_i)_{i\geq 1} \in \{k_1\cdots k_{m-1}(k_m+1), k_1\cdots k_m 0 \}$. All of this shows that $ (d_i)_{i\geq 1}$ satisfies \ref{(d)} (respectively \ref{(e)}) if and only if $\phi((d_i)_{i\geq 1}) $ satisfies \ref{(d)} (respectively \ref{(e)}). By induction, if $\phi^{n-1} ((d_i)_{i \geq 1})$ is the characteristic sequence associated with $x$, then $ \phi^n((d_i)_{i\geq 1}) = \phi(\phi^{n-1}((d_i)_{i\geq 1}))$ the characteristic sequence associated with $y$, for some $ y<-1$. This completes the proof. % periodic $\beta$-expansion of $l_{\beta}$ if and only if $ \phi((d_i)_{i\geq 1})$ is a periodic corrected  $x$-expansion of $l_x$, for some $ x $.  

%Now, suppose that $ (d_i)_{i \geq 1} = d(l_{\beta}, \beta)$. If $ (d_i)_{i \geq 1}$ is not periodic, $\phi(d(l_{\beta}, \beta)) $, 
\end{preu}

\subsection{Exchange of intervals}

This section establishes the exchange of the intervals $ ] -\infty, -\gamma_0 ] $ and $ ] -\gamma_n, -\gamma_{n+1} ]$ through the map $ \Upsilon$. 
We have the following results:

\begin{theo}\label{Th tau}
Let $n\in\N$. For any $\beta $ in $ ]-\gamma_n, -\gamma_{n+1} ]$, there is a unique $x \in ] -\infty, - \gamma_0 ]$ such that
 \begin{equation}
 d^{*}(l_{\beta}, \beta) = \phi^{n+1}(d^{*}(l_x, x)).
 \end{equation}
 Moreover, if $ \C_x$ codes the $x$-shift, then the support of the measure of maximal entropy on the $\beta$-shift is coded by $\phi^{n+1}(\C_x)$. %has the same language as $ \phi^{n+1}(S_x)$, where $ S_x $ denotes the $x$-shift. That is
 %\begin{equation}
  %A\in W \Leftrightarrow  A = \sigma^k(B), \text{ $ k\in\N$ and $ B\in \phi^{n+1}(S_x) $ }.
 %\end{equation}
 Let $ S_x^r$ and $ S_{\beta}^r$ denote the right one-sided subshifts associated with $ S_x$ and $ S_{\beta}$. We set
 \begin{equation*}
  t_n = \prod\limits_{k=-1}^{n-1}(1+\frac{1}{\beta^{l(u_k)}})-\frac{\beta^{l(u_n)}-2}{\beta^{l(u_n)-1}(\beta-1)}.
 \end{equation*}
Then, $ \phi^{n+1}(S_x^r)$ is the subset of $ S_{\beta}$ of $\beta$-representations consisting of real numbers of the interval $ [l_{\beta}, t_n ] $. The support of the measure of maximal entropy of the $\beta$-transformation is $ \underset{k \in \N}{\bigcup} T_{\beta}^k \left( [l_{\beta}, t_n ] \right)$. 
\end{theo}

\begin{theo}\label{Th surp}
Let $\Upsilon $ be the map from $ ] -\gamma_0, -1 [ $ to $] -\infty, -\gamma_0 ] $ such that 
\begin{equation}
d^{*}(l_x, x) = \phi^{n+1} \left( d^{*} (l_{\Upsilon(x)}, \Upsilon(x) )\right) \text{ if $ x \in ]-\gamma_n, -\gamma_{n+1} ] $ }.
\end{equation}
The map $\Upsilon $ is onto on its domain. Its restriction from $ ]-\gamma_n, -\gamma_{n+1} ] $ into $ ]-\infty, -\gamma_0 ] $ is one to one and continuous. 
\end{theo}

In the light of Theorem \ref{Th tau}, the map $ \phi^{n+1} $ transforms the $ x $-expansions of real numbers of the interval $ [l_x, r_x )$ into $\beta$-expansions of real numbers of $ [l_{\beta}, t_n )$.

\begin{theo}\label{th3}
Let $\beta <-1$. The $\beta$-shift is coded if and only if $ \beta \leq -\frac{1+\sqrt{5}}{2} $ and the $\beta$-expansion of the left endpoint $ l_{\beta}$ of the domain of $ T_{\beta}$ is not periodic with odd period. In all cases, the measure of maximal entropy $\mu_{\beta}$ is carried by a subset coded by a positive recurrent prefix language $P_{\beta}$.
It is the Champernowne measure of this prefix positive recurrent code $P_{\beta}$ such that for any $x \in P_{\beta}$, we have:
 \begin{equation}
  \mu_{\beta}([x]) = \dfrac{1}{|\beta|^{l(x)} \sum\limits_{y\in P_{\beta}}\frac{l(y)}{|\beta|^{l(y)}}}.\label{measCyl}
 \end{equation}
 \end{theo}
 In \cite{MR2974214}, the authors proved that $ (I_{\beta}, T_{\beta})$, is intrinsically ergodic. Both dynamical systems $ (I_{\beta}, T_{\beta})$ and right one-sided negative $\beta$-shift $ S_{\beta}^r $ endowed with the shift $ \sigma$, are conjugated. Thus, $ (S_{\beta}^r, \sigma) $ is intrinsically ergodic too. The following theorem extends this result to the two-sided negative $\beta$-shift endowed with the shift. 
 
 \begin{theo}\label{Theo 4}
 Let $ \beta < -1 $. Then, $ S_{\beta} $ is intrinsically ergodic. 
 \end{theo}
 
\begin{figure}[h!]
    \centering
\includegraphics[width=\textwidth]{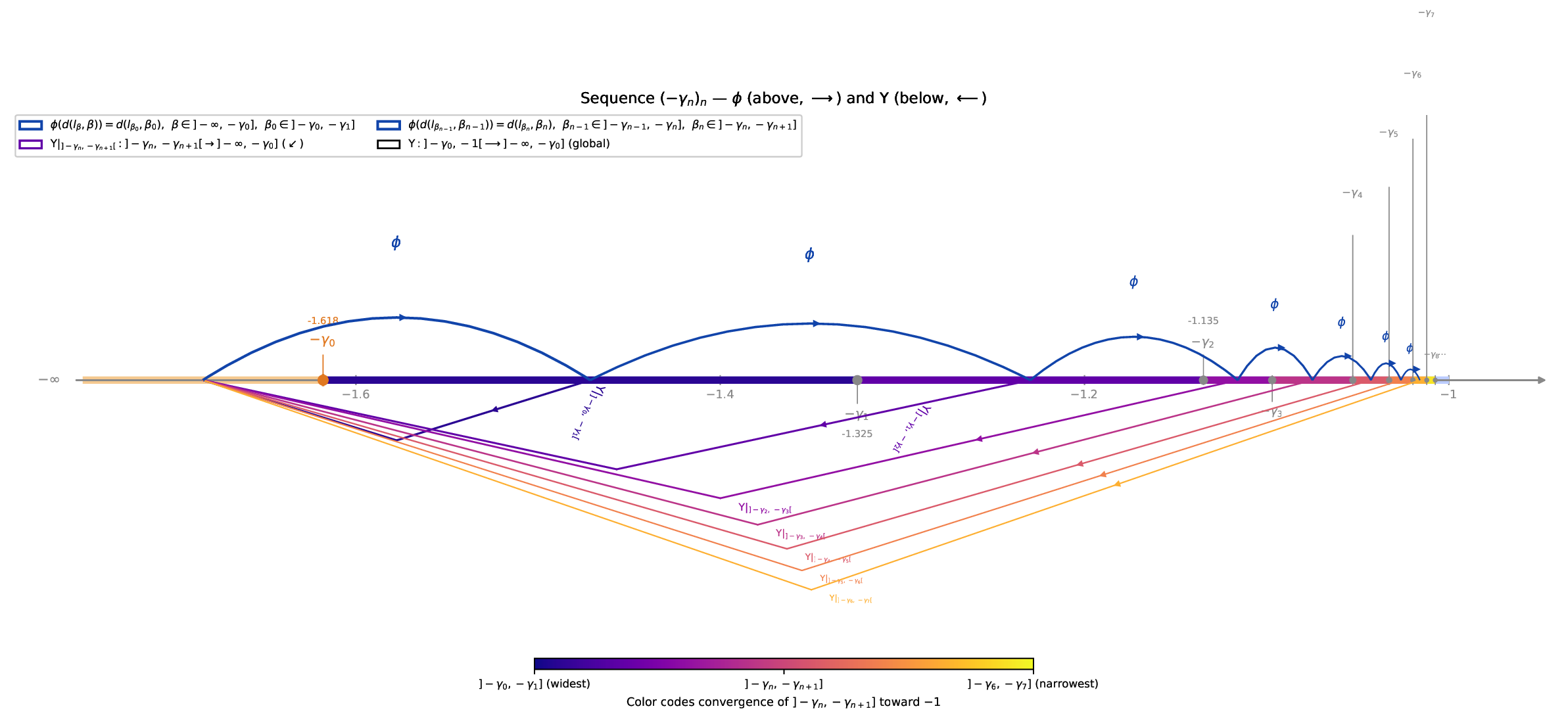}
\caption{Sequence of $\gamma_n$ and exchange of intervals.}
\end{figure}

In the rest of this paper, when $d(l_{\beta}, \beta)$ is periodic with odd period, we will use,   $S_{\beta}$ instead of $\tilde{S}_{\beta}$ as $\beta$-shift.

\section{Positive recurrent prefix code and negligible  cylinders}

Let $ (X, T) $ be a topological dynamical system. The variational principle describes the relationship between topological entropy and Kolmogorov entropy of a measurable dynamical system. We denote by $ h_{top}(T) $ the topological entropy. If $ \mu $ is a $ T$-invariant measure of $ X $, we denote by $ h_{\mu}(T) $ the usual metric entropy of $T$. In view of the variational principle, $ h_{top} $ coincides with the supremum of the metric entropy taken over all $T$-invariant measures. From \cite{MR0352411}, the $ \sup $ can be considered just over the set of ergodic measures.
\begin{equation*}
 h_{top}(T) = sup \{h_m(T) \vert m \in \M_T(X) \text{ is ergodic } \}.
\end{equation*}
In symbolic dynamics, setting $ T= \sigma$ (the shift on sequences), the $sup$ exists because $\sigma$ is an expansive map (see more details in \cite{MR0457675}). That is, we can find a real $ \theta > 0 $ such that for all $x$ and $y$ in any subshift $X$, $x\neq y$, there exists an integer $n $ which satisfies:
\begin{equation*}
 d( \sigma^{n}(x), \sigma^{n}(y)) \geq \theta 
\end{equation*}
where $d$ is the metric defined in \eqref{diadic 1} and \eqref{diadic 2}. It suffices to set $ \theta = \frac{1}{2}$.

Roughly speaking, a system is intrinsically ergodic if it has a unique measure of maximal entropy. The ergodicity of a measure implies the transitivity of its support. It is well known that coded systems are transitive. Bearing this in mind, and before adapting the study on coded systems done in \cite{MR858689} and \cite{MR939059} to the negative beta-shift, we begin by describing some codes we need to construct admissible sequences.

\subsection{Different codes to describe admissible sequences}
The set of factors of the $\beta$-shift is recognized by the automata given in Figure \ref{fig1} for the parameter $ \beta >1$, Figure \ref{fig2} when $\beta <- \frac{1+\sqrt{5}}{2} $ and Figure \ref{fig3} if $ \beta = - \frac{1+\sqrt{5}}{2} $. We denote by $ I_i $ the interval of integers such that $ I_{2i-1} = [0; d_{2i-1}-1]\cap \N $, $ I_{2i} = [ d_{2i}+1; d_1-1 ]\cap \N $ and $ D_i = [0 ; d_i-1 ]\cap \N $.

We begin this subsection by providing additional details (given in \cite{NguemaNdong2019}) on the construction of the code. Let $ (r_i)_{i \geq 1 } = \lim\limits_{x \rightarrow r_{\beta}}d(x, \beta) $.

\subsubsection{Case $ \beta > 1 $}
\begin{figure}[!ht]
      \centering
        \unitlength=4pt
  \[\atm\xymatrix{
   0\ar[r]_{a_1}\ar@(dl,dr)_{D_1}
             & 1\ar[r]_{a_2} \ar@/_0.4cm/[l]_{D_2}
             & 2\ar[r]_{a_3} \ar@/_0.8cm/[ll]_{D_3}
             & 3\ar[r]_{a_4} \ar@/_1.2cm/[lll]_{D_4}
             & 4\ar@/_1.6cm/[llll]_{D_5} \ar[r]_{a_5}
             & 5\ar@/_2cm/[lllll]_{D_6} \ar[r]_{a_6}
             & *{}
   }\]
    \caption{Automaton with $\beta> 1$.}\label{fig1}
    \end{figure}
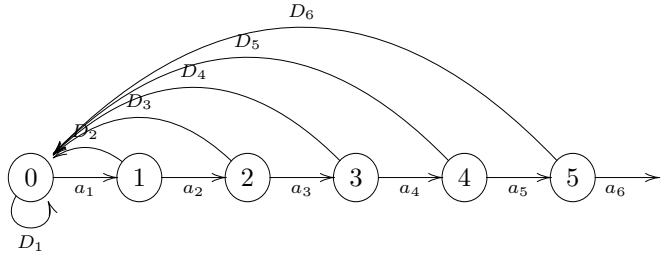

If $ \beta> 1 $, $(r_{i})_{i \geq 1}$ is the $ \beta$-expansion of 1. In this case, the $\beta$-shift is coded by the language
$ \{ r_1 r_2 \cdots r_{n}j \hspace{0.1cm}| \hspace{0.1cm} 0 \leq j <r_{n+1 }, \text{ $ n \in \N $ }\}$
with $ r_1\cdots r_n = \epsilon $ (the empty word) if $ n = 0 $. In this case, the set of factors recognizes the automaton of Figure \ref{fig1}.

\subsubsection{Case $ \beta \leq -\gamma_0$}

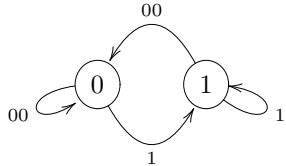
\begin{figure}[!ht]
      \centering
        \unitlength=4pt
  \[\atm\xymatrix{
   0\ar@(l,dl)_{00}\ar@/_0.8cm/[r]_{1}
             & 1\ar@(dr,r)_{1}\ar@/_0.8cm/[l]_{00}
             & *{}
   }\]
    \caption{Automaton with $ \beta = -\frac{1+\sqrt{5}}{2}$.}\label{fig3}
    \end{figure}

If $ \beta = -\frac{1+\sqrt{5}}{2}$, $ d(l_{\beta}, \beta) = \bullet 1(0)^{\infty} = \bullet 1\overline{0} $, the negative $\beta$-shift is coded by the prefix code $ \{1, 00\} $. The situation is summarized in Figure \ref{fig3}.

\begin{figure}[!ht]
      \centering
        \unitlength=4pt
  \[\atm\xymatrix{
   0\ar[r]_{d_1}\ar@(l,dl)_{I_1}
             & 1\ar@(dl,dr)_{d_1} \ar[r]_{d_2} \ar@/_0.4cm/[l]_{I_2}
             & 2\ar[r]_{d_3} \ar@/_0.8cm/[ll]_{I_3}
             & 3\ar[r]_{d_4} \ar@/^0.8cm/[ll]^{d_1} \ar@/_1.2cm/[lll]_{I_4}
             & 4\ar@/_1.6cm/[llll]_{I_5} \ar[r]^{d_5}
             & 5\ar@/^1.6cm/[llll]^{d_1} \ar@/_2cm/[lllll]_{I_6} \ar[r]^{d_6}
             & *{}
   }\]
    \caption{Automaton with $ \beta< -1$}. \label{fig2}
    \end{figure}
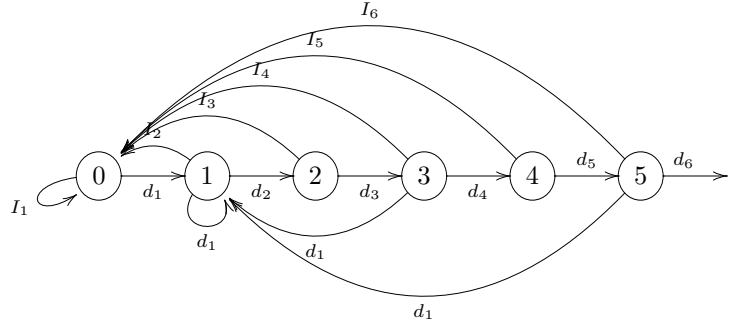

Suppose $ \beta < -1 $ and let $ d^{*}(l_{\beta}, \beta) = (d_i)_{i \geq 1} $. In \cite{NguemaNdong2019}, it is proved that $ \tilde{S}_{\beta}$ is coded if and only if  $ -\frac{1+\sqrt{5}}{2} \geq \beta $ and the $\beta$-expansion of $ l_{\beta} $ is not periodic with odd period (consequently,  $ S_{\beta}$ is coded if and only if $ -\frac{1+\sqrt{5}}{2} \geq \beta $). The situation is illustrated in Figure \ref{fig2}.

%In what follows, when $ d(l_{\beta}, \beta) $ is periodic with odd period, we will consider as negative $\beta$-shift, the set $ \tilde{S}_{\beta}$.
In light of Figure \ref{fig2}, if $ \beta < -\frac{1+\sqrt{5}}{2} $, $S_{\beta}$ is coded by the language $ \C_{\beta} $ consisting of words ending with a block of the form $ d_1 \cdots d_k j $, where $(-1)^{k+1}(d_{k+1} - j)< 0$. It is the set of words extendable to the right by any admissible sequence. These words are of the form $ x y $, where $ x $ is either the empty word or a concatenation of odd length initial segments of $(d_i)_{i \geq 1}$, that is, concatenation of words of the form $ d_1 \cdots d_{2p+1}$) and $ y =  d_1 \cdots d_k j $, where $(-1)^{k+1}(d_{k+1} - j)< 0$, $j\neq d_1$. We set
 
 \begin{equation*}
D= \{  d_1 \cdots d_{2k+1}, k \in \N \}.
\end{equation*}
Then, the concatenations of words of $ D $ are factors of words of the code. Note that all concatenations of $ D $ are admissible only if for any integer $ n \in \N^{\times} $, $ d_{2n} < d_1 $. It corresponds to the simplest case. More generally, we will put  $(d_n)_{n \geq 1}$ into a suitable form by introducing two sequences: one of odd positive integers $(2n_i-1)_{i \geq 1}$, and the other $ (p_i)_{i\geq 1}$ verifying $ 2n_i-1 > p_i \geq 1 $, $ d_{n_i-1+m} = d_m $ and $ (-1)^{2n_i+p_i}(d_{2n_i+p_i}-d_{p_i+1})< 0$. The blocks $ d_1 \cdots d_{p_i} $ are inserted after odd positions $ 2n_i-1 $. The sub-sequence $ (d_n)_{n \geq 2n_i} $ deviates from the sequence $(d_n)_{n \geq 1}$ at the index $ 2n_i+p_i$ for the first time since the index $ 2n_i$.

\begin{equation}
 d_1d_2\cdots = d_1 d_2 \cdots d_{2n_1-1}d_1 \cdots d_{p_1} d_{2n_1 + p_1} \cdots d_{2n_2-1}d_1 \cdots d_{p_2} d_{2n_2 + p_2} \cdots.
 \label{(6)}
\end{equation}

We adopt the convention that when $ n = 0$, $ d_1 \cdots d_n $ is empty. We now recall the construction of the codes introduced in \cite{NguemaNdong2019}. We set:

\begin{equation*}
\Gamma_0 = \{ d_1\cdots d_nj = y \in L_{\beta} | \text{ $ n\in \N $, $ d_1 \cdots d_{n+1-i} \prec \sigma^{i}(y)$, $ \forall i$ , $ 0 \leq i \leq n $}\};
\end{equation*}
\begin{equation}
\Delta_{0}^0 = \{ d_1 \cdots d_{2k-1} \vert 2n_i+p_i \leq 2k-1< 2n_{i+1}-1 \text{ and } i \in \N\}\label{Dood0}
\end{equation}
with $ n_0 = p_0 = 0 $; 
\begin{equation*}
E=\{ d_1\cdots d_{2n_i-1}; \text{ $ i \in \N^{\times}$ }\};
\end{equation*}
\begin{equation}
\C_{\beta}=\{xy \in L_{\beta} | x \in \{\varepsilon\}\cup D^*,\hspace{2mm} y \in \Gamma_0 \} \label{C};
\end{equation}
\begin{equation}
\Delta_{0,\beta} = \{ x y \in L_{\beta} | x \in \{\varepsilon\}\cup E^*, \hspace{2mm} y \in \Delta_0^0 \}\label{Dood}.
\end{equation}
That is, $\C_{\beta}$ and $\Delta_{0, \beta}$ are both codes satisfying the following properties:
\begin{itemize}
 \item First, for any $ x $ in $ \C_{\beta}$
 \begin{equation}
 d_1\cdots d_{l(x)-i} \prec\sigma^i (x),\label{property C}
\end{equation}
The fact that $ x $ can be extended to the right by any admissible word means that both $ x d_1 d_2 \cdots $ and $ x0d_1d_2 \cdots$ are admissible;
\item for any $ y $ in $\Delta_{0, \beta}$, both $ yd_1d_2\cdots $ and $ yd_1d_1d_2 \cdots$ are admissible. That is, each word of $\Delta_{0, \beta}$ can be extended to the right by any admissible sequence beginning with $ d_1$. %In the meaning of the alternating order, se
\end{itemize}
%and
%\begin{equation}
%d_1 \cdots d_{l(y)} \prec \sigma_D^i (y), \text{ $ \forall y \in \Delta_{0, \beta}$}.\label{property D}
%\end{equation}

We set (see more details in \cite{NguemaNdong2019}) $ B_i = d_1d_2\cdots d_{2n_i-1} $ and
\begin{equation}
J(0) = \{ t \vert p_t < 2n_1-1 \}, \label{J0}
\end{equation}
and for all $ i \in \N^{\times}$,
\begin{equation}
J(i) = \{ t \vert 2n_i-1 \leq p_t < 2n_{i+1}-1 \} \label{Ji};
\end{equation}
and we denote by $ \Delta_{i,\beta} $ the set of words $ x $ such that
\begin{equation}
\begin{cases}
                    x = B_{t_1}\cdots B_{t_m}, \\
                    p_{t_k} \leq 2n_{t_{k+1}}-1, \\
                     t_m \in J(i), \\
                      t_k \not\in J(i) \\
                    \text{for $ k\neq m $ }, p_{t_m}< 2n_{t_1}-1.
                    \end{cases} \label{Di}
\end{equation}
The language $ \C_{\beta} $ is a prefix code.

\begin{exple}\label{ex1}
Let $ \beta<0 $ be the algebraic integer associated with the polynomial
\begin{equation}
P(X)= X^{15}+3X^{14}-X^{12}+2X^{11}-X^{6}-X^{5}+2X^{4}-X^{3}-2X^{2}+2X+1.
\end{equation}
Then, $d(l_{\beta}, \beta) = \bullet 2012121201200 \overline{21} $, \\
$ \Delta_{1, \beta} = \{201, \overline{2012121}^k 20121, \overline{2012121}^k2012121201200\overline{21}^p| p, k\geq 0 \} $, \\
 $ \Delta_{2, \beta} = \{2012121 \} $, $ \Delta_{0,\beta} = \{2\}$, $ \Gamma_0 = \{0, 1, 21, 200 \} $, \\$ \Gamma_1 = \{ x 200 \vert x \in \Delta_{1, \beta}^* \} $ and $ \C_{\beta} = \{ x y\vert x \in \Delta_{0,\beta}^*, y \in \Gamma, l(y) \geq 2 \} \cup \Gamma $.
\end{exple}

\begin{exple}
Let $ \beta < -1 $ be the algebraic integer of the polynomial
\begin{equation}
P(X) = X^{14}+2X^{13}-2X^{12}+X^{11}+X^{10}-X^{9}+X^{8}-X^{7}+X^{6}-2X^{5}+X^{4}+X^{3}-2X^{2}+1.
\end{equation}
Then, $ d(l_{\beta},\beta) = \bullet 2012121201200 \overline{1} $. \\
 $ \Delta_{1, \beta} = \{ 201, \overline{2012121}^k 20121 | k \geq 0 \} $, $\Delta_{2, \beta} = \{2012121 \} $;\\
 $ \Delta_0^0 = \{2, 2012121201200(11)^k| k \geq 0\} $\\ $ \Delta_0^1 = \{xy2012121201200(11)^k| k \geq 0, x \in \Delta_{1, \beta}^{*}, y \in \Delta_{2, \beta}^* \} $,\\
$ \Gamma_0 = \{0,1, 21, 200, 2012121201200(11)^k10| k \geq 0 \} $; \\
$ \Gamma_1 = \{ x200, xy 2012121201200(11)^k10| k \geq 0, x \in \Delta_{1, \beta}^*, y \in \Delta_{2,\beta}^* \} $.
\end{exple}

%If $ -\frac{1+\sqrt{5}}{2}< \beta<-1 $, $S_{\beta}$ is not a coded system and one has $ \C_{\beta} = \{0\}$. In this case, Figure \ref{fig1} does not suffice for the study the factors of $ S_{\beta} $. So, the question to know how to generate concatenations of words of $D$ seems necessary.

\subsubsection{Case $-\gamma_0 < \beta < -1$}

If $ -\frac{1+\sqrt{5}}{2}< \beta<-1 $, $S_{\beta}$ is not a coded system and $ \C_{\beta} = \{0\}$. In this case, Figure \ref{fig1} does not suffice for the study of factors of $ S_{\beta} $. So, the question to know how to generate concatenations of words of $D$ seems necessary.

 \begin{prop}\label{ppp1}
  Let $ \beta <- 1 $ and $ d^{*}(l_\beta, \beta) = (d_i)_{i \geq 1}$. We denote by $ D $ the subset of $ S_{\beta} $ consisting of admissible concatenations of words of the type $ d_1 \cdots d_{2n-1} $ and $m$ an ergodic measure on $ S_{\beta} $ of maximal entropy. If the measure $m$ is carried by $D$, then: 
 \begin{equation}
 h_m(S_{\beta}) \leq \log \frac{1+\sqrt{5}}{2}.
 \label{(26)}
 \end{equation}
 \end{prop}

\begin{preu}
Let $F$ be the set of words on $ \{0, 1, \cdots, d_1\} $ generated by the initial segments of $ (d_i)_{i \geq 1}$ with odd lengths and $ F = F^{\times} \cup \{\varepsilon\} $. Then, $ D \subset F $. But

\begin{equation}
F^{\times} = d_1F \cup  (d_1d_2d_3F) \cup (d_1d_2d_3d_4d_5F) \cup \cdots.
 \end{equation}
 This implies that the number $f_n$ of words with length $n$ of $F$ satisfies:
\begin{align*}
  f_n &= f_{n-1} + f_{n-3} + f_{n-5} + \cdots. \\
      &= f_{n-1} + f_{n-2}.
 \end{align*}
And thus, $ \frac{1}{n}\log f_n $ tends to $ \log \frac{1+\sqrt{5}}{2} $. Since $ D \subset F $, we obtain the result.
\begin{equation}
 h_m(S_{\beta}) \leq \log \frac{1+\sqrt{5}}{2}.
\end{equation}
\end{preu}

\begin{prop}\label{pp2}
Let $ \beta \in [ -\gamma_{n}, -\gamma_{n+1} )$. Then $ d^{*}(l_\beta, \beta) $ is the image under $ \phi^n$ of the characteristic sequence associated with $x$, for some $ x $ satisfying $ -\gamma_1 > x \geq -\gamma_0$.
\end{prop}

We recall that $ d(l_{-\gamma_n}, -\gamma_n)= d^{*}(l_{-\gamma_n}, -\gamma_n) = u_{n}(u_{n-1}u_{n-1})^{\infty}$.

\begin{preu}
 Consider $ \beta \in [ -\gamma_{n}, -\gamma_{n+1} ) $. Then:
\begin{equation*}
d(l_{-\gamma_n},-\gamma_n) \preceq d^{*}(l_\beta,\beta) \prec d(l_{-\gamma_{n+1}}, -\gamma_{n+1}). 
\end{equation*}
So, there exists $ n_1 $ such that $ d(l_\beta, \beta) $ starts with $ u_n(u_{n-1})^{2n_1}u_n $. The word $ (u_{n-1})^{2n_1} $ is the longest concatenation of $ u_{n-1}$ which follows $ u_n $ in an admissible sequence. But all sequences of length $ l( u_n(u_{n-1})^{2n_1}u_n )=2l(u_n)+2n_1l(u_{n-1}) $ are greater than $ u_n(u_{n-1})^{2n_1}u_n $ (in the meaning of the alternating order) and after $ u_n $ one must have an even number of $ u_{n-1} $. Hence, there exists a bounded sequence $ (n_i)_{i \geq 1} $, $ n_i \leq n_1 $ such that  
 \begin{equation*}
  d^{*}(l_{\beta}, \beta) = u_n(u_{n-1})^{2n_1}u_n(u_{n-1})^{2n_2}u_n (u_{n-1})^{2n_3} \cdots = \phi^n(1(0)^{2n_1}1(0)^{2n_2} \cdots).
 \end{equation*}
 
 Let us show that $ 1(0)^{2n_1}1(0)^{2n_2} \cdots$ is the characteristic sequence associated with $x$. Then for any $ k \in \N$,
\begin{equation}
 1(0)^{2n_1}1(0)^{2n_2} \cdots \preceq (0)^k1(0)^{2n_i}1(0)^{2n_{i+1}} \cdots. \label{(A)}
\end{equation}
 Otherwise, since $\phi$ is increasing, there would exist $k $ such that
 \begin{equation*} 
 u_{n-1}^ku_n(u_{n-1})^{2n_i}u_n(u_{n-1})^{n_{i+1}}\cdots \preceq  u_n(u_{n-1})^{2n_1}u_n(u_{n-1})^{2n_2}u_n (u_{n-1})^{2n_3} \cdots .
 \end{equation*} 
 This is absurd because $ u_n(u_{n-1})^{2n_1}u_n(u_{n-1})^{2n_2}u_n (u_{n-1})^{2n_3} \cdots $ is the characteristic sequence associated with $\beta$.
 
 Moreover, 
 \begin{equation}
 1(00)^{n_1}1(00)^{n_2}\cdots \prec u_n(u_{n-1})^{2n_1}u_n(u_{n-1})^{2n_2} \cdots \prec \phi^{\infty}(1) .\label{(B)}
 \end{equation}
 From \ref{(A)} and \ref{(B)}, there is $x$ such that $f_{x}( 1(00)^{n_1}1(00)^{n_2}\cdots) = \frac{x}{1-x}$. 
 
 In the proof of Corollary \ref{Cor Tau}, we have seen that $ 1(00)^{n_1}1(00)^{n_2}\cdots $ satisfies the \ref{(d)} (respectively \ref{(e)}) if and only if the same holds for $d^{*}(l_{\beta}, \beta)$.
 Hence the result follows.
\end{preu}

As a consequence of Proposition \ref{pp2}, we have the following corollary:

\begin{coro}\label{copp2}
Let $\beta $ be a real number, with $\beta \in ]-\gamma_{n}, -\gamma_{n+1}]$. Then  there is a unique $ x \in ]-\infty, -\gamma_0]$ such that $d^{*}(l_{\beta}, \beta) = \phi^{n+1}(d^{*}(l_x, x))$.
\end{coro}

\begin{preu}\textbf{Corollary \ref{copp2}}
\begin{itemize}
 \item[(1)] Let $ \beta \in ] -\gamma_0, -\gamma_1 ]$. We know that $ d(l_{-\gamma_0}, -\gamma_0) = \bullet 1(0)^{\infty} $ and $ d(l_{-\gamma_1}, -\gamma_1) = \bullet 100(11)^{\infty} $. Thus, for some sequences of integers $(s_i)_{i\geq 1}$ and $(q_i)_{i\geq 1}$,
\begin{equation}
\begin{aligned}
d^{*}(l_{\beta}, \beta) &= 1001(11)^{s_1}100(11)^{s_2}100(11)^{s_3}100(11)^{s_4}1\cdots \\
                    &= \phi \left( 10(00)^{s_1}1(00)^{s_2}1(00)^{s_3}1(00)^{s_4}\cdots \right).
\end{aligned}
\end{equation}
or
\begin{equation}
\begin{aligned}
d^{*}(l_{\beta}, \beta) &= 1(00)^{q_1+2}1(00)^{q_2}1(00)^{q_3}1(00)^{q_4}\cdots \\
                    &= \phi\left((q_1+2)q_2q_3q_4\cdots\right).
\end{aligned}
\end{equation}
From Corollary \ref{Cor Tau}, $10(00)^{s_1}1(00)^{s_2}1(00)^{s_3}1(00)^{s_4}\cdots$ (which is smaller than $d(l_{-\gamma_0}, -\gamma_0) = 1(0)^{\infty}$ in the meaning of the alternating order) is the characteristic sequence associated with $x$, for some $ x $. The same is true for $(q_1+2)q_2q_3q_4\cdots$. In all cases, $x \geq -\gamma_0$.
\item[(2)] For $\beta \in ]-\gamma_{n}, -\gamma_{n+1}]$, the result holds using Proposition \ref{pp2} and (1) of the present proof.
\end{itemize}
\end{preu}

\begin{rem}\label{remg}
 Let $\beta\in ]-\gamma_0, -\gamma_1]$. We know that we can find $x\in ]-\infty, -\gamma_0]$ such that $d^{*}(l_{\beta}, \beta) = \phi( d^{*}(l_x, x))$. It is easy to check that $\Delta_{0, \beta} = \phi(\C_x)$ where $ \Delta_{0,\beta}$ is defined in \eqref{Dood} and $\C_x$ the code of $S_x$ given in \eqref{C}. More generally, if $\beta\in ]-\gamma_{n-1}, -\gamma_{n}]$, $\Delta_{n-1, \beta} = \phi^{n}(\C_{x})$, with $\phi^{n}(d^{*}(l_x, x)) = d^{*}(l_{\beta}, \beta)$. 
\end{rem}

\begin{exple}
Let $d^{*}(l_x, x) = 2\overline{1}$.
\begin{equation*}
\begin{aligned}
&\Delta_{0,x} = \{2(11)^k; k\in \N\};\\
&\C_{x} = \{0, 1, 2(11)^k10, A2(11)^k10| A\in \Delta^*_{0, x }\} = \{0, 1,  A10 | A \in \Delta^*_{0,x} \}; \\
&d^{*}(l_{\beta}, \beta)  = \phi(d^{*}(l_x, x)) = 10000\overline{100};\\ 
 &\Delta_{1, \beta} = \{ 10000(100100)^k; k \in \N \} = \phi(\Delta_{0, x});\\
& \Delta_{0, \beta} = \{1, 100, A1001; A \in \Delta_{1, \beta}^*\} = \phi(\C_x)  
\end{aligned}
\end{equation*}
\end{exple}

In the interval $ ] \gamma_1, \gamma_0 [$, there is an increasing sequence of real numbers $ (\alpha_n)_{n \geq 1}$ such that $ \alpha_n$ is the algebraic number of $ X^{2n+1}-X^{2n}-X^{2n-1}+1 $ and $ \phi(d^{*}(l_{-n}, -n)) = d^{*}(l_{-\alpha_n}, -\alpha_n) $. It is easy to see that $ \gamma_0$ is its limit when $ n $ tends to infinity. For any $ \beta \in [-\alpha_{n}, -\alpha_{n+1}[ $, the length of the longest string of zeros which follows 1 is $2n$.
%\begin{prop}

%\end{prop}

\begin{lem}\label{lm1}
Let $ \beta $ be a real number such that  $ -1> \beta \geq -\gamma_0 $. Then  $ Card (\Delta_{0,\beta} ) \geq 2$ if and only if $ \beta <-\gamma_1 $, where $ Card( \Delta_{0,\beta} ) $ is the cardinality of $ \Delta_{0,\beta} $.
\end{lem} 

\begin{preu}
\begin{itemize}
\item Suppose $ Card( \Delta_{0,\beta} ) \geq 2 $. Note that when $ -\gamma_0$  is smaller than $ \beta $, the length of any string of zeros which appears in the characteristic sequence associated with $\beta$ is even. Moreover, remark that if $ \Delta_{0,\beta} $ contains a word of length 3 (that is 100), then 10000 is admissible and thus $ d^{*}(l_{\beta}, \beta)\prec d(l_{-\gamma_1}, -\gamma_1)$. Indeed $ d(l_{-\gamma_1}, -\gamma_1) = 100\overline{11} $ and $ 10000 \prec 10011$.

Now suppose that $ 100 \notin \Delta_{0,\beta} $. Then $ d^{*}(l_{\beta}, \beta) $ starts with $ 10011 $. For $ u \in \Delta_{0,\beta} $ with length at least $ 5 $, if $ 10011 \notin\Delta_{0,\beta}$, $ 10011u \in \Delta_{0,\beta} $ (by definition of $ \Delta_{0,\beta}$). The longest sequence of zeros is 00. So, $u$ ends with a sequence of the type $100(1)^t$ for some integer $ t $. Since $ 1\in \Delta_{0, \beta}$ (all concatenations of words of $ \Delta_{0,\beta}$ are admissible), it follows that $ 100(11)^{\infty}$ is admissible, that is $ d^{*}(l_{\beta}, \beta) \prec 100\overline{11} $. Thus, from Proposition 2 of \cite{NguemaNdong20161}, $ \beta <- \gamma_1 $.

\item Suppose $ -\gamma_1 > \beta \geq -\gamma_0 $.
\begin{equation}
1(0)^{\infty} \prec d^{*}(l_{\beta}, \beta) \prec 100(11)^{\infty}.
\end{equation}
Then, $ d^{*}(l_{\beta}, \beta) = 100(11)^{t_1}00(1)^{t_2}00(1)^{t_3} \cdots $. If $ t_1 = 0$, $100 \in \Delta_{0,\beta}$. If $t_1 \neq 0$, $ 100(11)^{t_1} \in \Delta_{0,\beta}$. Consequently, $Card(\Delta_{0,\beta}) \geq 2$.
\end{itemize}
\end{preu} 

\begin{lem}\label{lm2}
Let $ \beta $ be a real number such that $ -\gamma_1 > \beta \geq -\gamma_0$. Then, the topological entropy of the free monoid generated by $\Delta_{0,\beta}$ is larger than that of the free monoid generated by $ \Delta_{n,\beta}$, with $ n \geq 1$.
\end{lem}

\begin{preu}
From Lemma \ref{lm1}, $ Card(\Delta_{0,\beta}) \geq 2 $. That is, $ \Delta_{0,\beta} \neq \{ d_1 \} $. Remark that if $ d_1 \cdots d_{2n_i-1} \notin L_{\Delta_{0,\beta}^*}$, then for all $ t \geq i $, $ d_1 \cdots d_{2n_t-1} \notin L_{\Delta_{0,\beta}^*}$. If such an integer $ i $ is minimal, $ d(l_{\beta}, \beta)$ is an infinite concatenation of two consecutive words (with respect to the alternating order) $ U_0 = d_1 \cdots d_{2n_i-1}$ and $ V_0 = d_1 \cdots d_{2n_i-2}(d_{2n_i-1}-1)0 $ or $ V_0 = d_1\cdots d_{2n_i-3}(d_{2n_i-2}+1) $ (see the proof of Theorem 2 of \cite{NguemaNdong2019}). In an infinite admissible sequence, $ U_0 $ and $ V_0$ are followed by $ U_0 $ or $ V_0$. Thus, we can find in $\Delta_{0,\beta} $ a word $ x \neq d_1 $ with length smaller than $ l(U_0)$ and if $ \Delta_{i_0, \beta} $ is the language which contains $ U_0 $, then $ \Delta_{i_0, \beta}^* \subset \{U_0, V_0 \}^{*} $.

Let $ \log\beta_1 $ be the entropy of $ \{U_0, V_0\}^* $ endowed with the shift. We have
\begin{equation}
\begin{aligned}
1 &= \frac{1}{\beta_1^{l(U_0)}}+\frac{1}{\beta_1^{l(V_0)}}\\
  &= \sum\limits_{n \geq 0} \frac{1}{\beta_1^{nl(U_0)+l(V_0)}}\label{U0}
\end{aligned}
\end{equation}
The entropy of $ \Delta_{i_0, \beta}^*$ is smaller than $ \log\beta_1$ since $ \Delta_{i_0, \beta}^* \subset \{U_0, V_0 \}^{*} $.

We have seen that there is a word $ x $ in $\Delta_{0,\beta} $ such that $ x \neq d_1$ and $l(x)< l(U_0)$. Then $ \{d_1, x \}^* \subset \Delta_{0,\beta}^* $. Let $ \log \beta_2 $ be the entropy of $ \{d_1, x \}^*$. We have
\begin{equation}
1 = \frac{1}{\beta_2}+\frac{1}{\beta_2^{l(x)}}.
\end{equation}
Let $ d_{\beta_2}(1) $ be the $\beta_2$-expansion of 1. Then, $ d_{\beta_2}(1) = 1(0)^{l(x)-1}1 $. So, because $l(x)< l(U_0)$, it follows that $ (0)^{l(V_0)-1}\overline{1(0)^{l(U_0)-1}}$ is an infinite word of the $\beta_2$-shift. This implies that $ \sum\limits_{n \geq 0} \frac{1}{\beta_2^{nl(U_0)+l(V_0)}}<1 $. Consequently $\beta_1<\beta_2$ since $ \beta_1 $ is the largest real number satisfying \eqref{U0} and the map $ z\mapsto \sum\limits_{n \geq 0} \frac{1}{z^{nl(U_0)+l(V_0)}}$ on $ \R_{+}^{*}$ decreases. Hence, the result follows.
\end{preu}

\begin{rem}\label{nlem1}
If $ -\gamma_{n+1} > \beta \geq -\gamma_n $ then, for all $ k< n$,  $\Delta_{k, \beta} = \{ u_k \} $ (see the proof of Lemma 7 of \cite{NguemaNdong2019}).
\end{rem}

\subsubsection{Recurrent prefix codes}
 
\begin{lem}\label{lmm3}
Let $ \beta < -1 $. Then:
\begin{itemize}
\item if $ \beta < -\gamma_0 $, then $ \sum\limits_{x \in \C_{\beta}}\frac{l(x)}{|\beta|^{l(x)}} <+\infty $ and $ \sum\limits_{ x \in \C_{\beta}} \frac{1}{|\beta|^{l(x)}} = 1$;
\item if $ -\gamma_n \leq \beta < -\gamma_{n+1}$, then $ \sum\limits_{x \in \Delta_{n,\beta}}\frac{l(x)}{|\beta|^{l(x)}} < +\infty $ and $ \sum\limits_{ x \in \Delta_{n,\beta}} \frac{1}{|\beta|^{l(x)}} = 1$.
\end{itemize}
\end{lem}

\begin{rem}\label{remH}
The coefficients of the expansion in the formal power series of $ \dfrac{1}{\prod\limits_{k \geq 0}(1-\sum\limits_{x \in \Delta_{k, \beta}}z^{l(x)})} $ count the admissible concatenations of words  of the type $ d_1 \cdots d_{2k+1} $, with $k \in \N$. Thus, the smallest pole in modulus of $ \dfrac{1}{\prod\limits_{k \geq 0}(1-\sum\limits_{x \in \Delta_{k, \beta}}z^{l(x)})} $ is  greater than $\frac{2}{1+\sqrt{5}}=\frac{1}{\gamma_0}$. The equality holds only if all initial segments of odd lengths belong to $ \Delta_{0, \beta}$. That is  $\Delta_{0, \beta} = \{ d_1\cdots d_{2n+1}, n \in \N\}$ and $ \sum\limits_{x \in \Delta_{k, \beta}}z^{l(x)} = 0 $ if $ k \geq 1 $.
\end{rem}

\begin{preu 5}

If $ H_n $ denotes the number of words of length $ n $ in $ L_{\beta}$, it follows from Proposition 1 of \cite{NguemaNdong20161} that:
\begin{equation}
H_n = \sum\limits_{k=1}^{n}(-1)^k(d_{k-1}-d_k)H_{n-k}+1.
\end{equation}
Using Theorem 2 and Theorem 3 of \cite{NguemaNdong2019}, in the sense of formal power series, we have the following equation:
\begin{equation}
1-\sum\limits_{n \geq 1} (-1)^n(d_{n-1}-d_n)z^n = (1+z)(1-\sum\limits_{x\in \C_{\beta}}z^{l(x)})\prod\limits_{k \geq 0}(1-\sum\limits_{x \in \Delta_{k, \beta}}z^{l(x)}). \label{lab1}
\end{equation}
The left power series vanishes at $ -\frac{1}{\beta} = \frac{1}{|\beta|}$ which is its smallest root in modulus. That is $ \sum\limits_{n \geq 1}\frac{d_{n-1}-d_n}{\beta^n} = 1$ and so: 
%\begin{equation}
$ \sum\limits_{x\in \C_{\beta}}\frac{1}{|\beta|^{l(x)}} = 1 $
%\end{equation}
or there exists $ n \in \N $ such that $  \sum\limits_{x\in \Delta_{n, \beta}}\frac{1}{|\beta|^{l(x)}}=1 $.

From Proposition \ref{ppp1}:
\begin{itemize}
\item[(1)] If $ \beta <- \gamma_0 $, $S_{\beta}$ is coded by $ \C_{\beta} $. Moreover, the entropy of the system generated by the language $ \{ d_1\cdots d_{2n+1}| n \in \N\} $ is $\log \gamma_0 $. Thus the subsystem of admissible sequences which are concatenations of words of the type $ d_1\cdots d_{2n+1}$ has entropy smaller than $ \log \gamma_0 $. Then
\begin{equation}
 \prod\limits_{k \geq 0}(1-\sum\limits_{x \in \Delta_{k, \beta}}\frac{1}{|\beta|^{l(x)}}) \neq 0 \text{ and } 1-\sum\limits_{x\in \C_{\beta}}\frac{1}{|\beta|^{l(x)}} = 0.
\end{equation}
\item[(2)] If $ -\gamma_n \leq \beta <- \gamma_{n+1} $, $\C_{\beta}=\{0\}$, for $ i < n$, $ \Delta_{i, \beta} = \{u_i\} $ and for $ i > n $, $ \Delta_{i,\beta}^* \subset L_{\Delta_{n,\beta}^*}$ (see the proof of Lemma 7 of \cite{NguemaNdong2019}). Thus, the entropy of $ \Delta_{n,\beta}^* $ is greater than that of the subsystem of  concatenations of words of the sets $ \Delta_{i,\beta} $ with $ i \geq n+1 $. The coefficients of the expansion in the formal power series of $ \dfrac{1}{\prod\limits_{k \geq n+1}(1-\sum\limits_{x \in \Delta_{k,\beta}}z^{l(x)})} $ count admissible concatenations of words of the sets $ \Delta_{i,\beta} $ with $ i \geq n+1 $. Then $ 1-\sum\limits_{x \in \Delta_{n, \beta}}\frac{1}{|\beta|^{l(x)}} = 0$ and  $ \prod\limits_{k \geq n+1}(1-\sum\limits_{x \in \Delta_{k,\beta}}\frac{1}{|\beta|^{l(x)}}) \neq 0 $.
\end{itemize}
Now, we have
\begin{equation}
\begin{aligned}
\sum\limits_{x \in \C_{\beta}}\frac{1}{|\beta|^{l(x)}} = 1 &\text{ if $ \beta <- \gamma_0 $}\\
\sum\limits_{x \in \Delta_{n,\beta} }\frac{1}{|\beta|^{l(x)}} = 1 &\text{ if $ -\gamma_n \leq \beta < -\gamma_{n+1} $ }
\end{aligned}\label{eqCD2}
\end{equation}
Using the derivatives of the formal powers series in \eqref{lab1} and the relation \eqref{eqCD2}, it follows that:
\begin{equation}
\sum\limits_{n \geq 1}n\dfrac{(d_{n-1}-d_n)}{\beta^n} = 
\begin{cases}
(1+\frac{1}{|\beta|}) \sum\limits_{x \in \C_{\beta}}\frac{l(x)}{|\beta|^{l(x)}} \prod\limits_{k \geq 0} \left(1-\sum\limits_{x \in \Delta_{k,\beta}}\frac{1}{|\beta|^{l(x)}}\right) &\text{ if $ \beta <- \gamma_0 $ } \\
(1-\frac{1}{\beta^2})\sum\limits_{x\in \Delta_{n,\beta}}\frac{l(x)}{|\beta|^{l(x)}}\prod\limits_{k \neq n} \left(1-\sum\limits_{x \in \Delta_{k,\beta}}\frac{1}{|\beta|^{l(x)}}\right) &\text{ if $ -\gamma_n \leq \beta < -\gamma_{n+1} $}
\end{cases} 
\end{equation}
Since $ (d_{n-1}-d_n)_{n \geq 1} $ is bounded, it follows that 
\begin{equation}
\begin{aligned}
\sum\limits_{x \in\C_{\beta}}\frac{l(x)}{|\beta|^{l(x)}} <+\infty &\text{ if $ \beta <- \gamma_0 $ } \\
\sum\limits_{x \in \Delta_{n,\beta} }\frac{l(x)}{|\beta|^{l(x)}} < +\infty &\text{ if $ -\gamma_n \leq \beta < -\gamma_{n+1} $ }.
\end{aligned}
\end{equation}
\end{preu 5}

\begin{lem}\label{Plem}
Let $ P_{\beta} $ be the language
\begin{equation}
P_{\beta} = \begin{cases}
    \C_{\beta} &\text{ if $ \beta < -\frac{1+\sqrt{5}}{2}$ }\\
     \Delta_{n,\beta} &\text{ if  $ -\gamma_n \leq \beta < - \gamma_{n+1} $ },
  \end{cases}
\end{equation}
$ \rho_{P_{\beta}} $ and $\rho_{P_{\beta}^*}$ denote radii of $ \sum\limits_{x\in P_{\beta}} z^{l(x)} $ and $ \sum\limits_{x\in P_{\beta}^*} z^{l(x)} $ respectively. Then $ \rho_{P_{\beta}^*} <\rho_{P_{\beta}}$.
\end{lem}
It is obvious that $\rho_{P_{\beta}^*} \leq \rho_{P_{\beta}}$. In the following proof of Lemma \ref{Plem}, we show that the equality cannot hold.

\begin{preu} \textbf{of Lemma \ref{Plem}}
\begin{itemize}
\item[•] Suppose $ \beta < -\frac{1+\sqrt{5}}{2}$. Let $c_n$ be the number of words of length $ n$ in $\C_{\beta}$. From \eqref{lab1}, 
\begin{equation*}
1-\sum\limits_{n \geq 1}c_nz^n = \dfrac{1+\sum\limits_{n\geq 1}(-1)^n(1+d_n)z^n}{\prod\limits_{k \geq 0}(1-\sum\limits_{x \in \Delta_{k,\beta}}z^{l(x)})}. 
\end{equation*}
The radius $ \rho_{\C_{\beta}}$ is the smallest pole in modulus of $1-\sum\limits_{n \geq 1}c_nz^n$, that is the smallest root in modulus of all power series $ 1-\sum\limits_{x \in \Delta_{k,\beta}}z^{l(x)}$. But if $ z_0$ is the root of $ 1-\sum\limits_{x \in \Delta_{k,\beta}}z^{l(x)} $ with the smallest modulus, then $-\log |z_0|$ is the entropy of the monoid generated by $ \Delta_{k,\beta} $. Since, for any $ k $, $\Delta_{k, \beta} \subset \{d_1\cdots d_{2n+1} | n \in \N \}$ and the entropy of the monoid generated by $ \{d_1\cdots d_{2n+1} | n \in \N \}^{\N}$ is equal to $\log \gamma_0$, it follows that:
\begin{equation*}
-\log \rho_{P_{\beta}} \leq \log \gamma_0 < \log |\beta|. 
\end{equation*} 
We have the result using the fact that $ \C_{\beta} $ codes the system and then, $S_{\beta}$ and $\C_{\beta}^*$ have the same entropy $ \log |\beta| = -\log \rho_{\C_{\beta}^*}$.
%and Remark \ref{remH}, and since $(d_n)_{n \geq 1}$ is bounded, the radius of the formal power series
% is larger than the $\frac{2}{1+\sqrt{5}}$. More precisely, the smallest zero of $\prod\limits_{k \geq 0}(1-\sum\limits_{x \in \Delta_{k,\beta}}z^{l(x)}) $ belongs to $ [\frac{2}{1+\sqrt{5}}, +\infty ) $. Moreover, if $\beta< -\gamma_0$, $S_{\beta}$ is coded by $\C_{\beta}$. Thus, the language of $\C_{\beta}^{*}$ equal to the language of $S_{\beta}$. So, $\rho_{P_{\beta}^*} = \frac{1}{|\beta|} < \frac{1}{\gamma_0}$.
%Thus $ \rho_{P_{\beta}^*} < \frac{2}{1+\sqrt{5}} \leq \rho_{P_{\beta}}$.
\item[•] If $ -\gamma_n \leq \beta < -\gamma_{n+1} $, $ \Delta_{i, \beta} = \{u_i\}$ if $ 0 \leq i < n $, $ \C_{\beta} = \{0\} $ and $ P_{\beta} = \Delta_{n, \beta} = \phi^{n+1}(\C_x)$, where $d(l_{\beta}, \beta) = \phi^{n+1}(d(l_x, x)$. We know that the topological entropy of the language generated by the free monoid $\C_x $ is strictly larger than the entropy of $ \Delta_{0, x}^* $. This implies that the entropy of $\Delta_{n, \beta}^* = \phi^{n+1}(\C_x)^*$ is strictly larger than that of $\Delta_{k, \beta}^*$ for any $k\geq n+1$. The radius $ \rho_{P_{\beta}} $ is the smallest pole in modulus of
\begin{equation*} 1-\sum\limits_{x\in P_{\beta}}z^{l(x)} = \dfrac{1+\sum\limits_{n\geq 1}(-1)^n(1+d_n)z^n}{(1-z)\prod\limits_{k = 0}^{n-1}(1-z^{l(u_k)})\prod\limits_{k \geq n+1}(1-\sum\limits_{x \in \Delta_{k,\beta}}z^{l(x)})}.
\end{equation*} 
That is, $ \rho_{P_{\beta}} = |z_0|$, where $ z_0 $ is the smallest root of the sums $ 1-\sum\limits_{x \in \Delta_{k,\beta}}z^{l(x)} $, $ k \geq n+1$. But, if so, $-\log |z_0| $ is equal to the entropy of $ \Delta_{k, \beta}^*$, which is strictly smaller than $-\log \rho_{P_{\beta}^*} $, entropy of $ P_{\beta}^* $. It follows that $ \rho_{P_{\beta}^*} < \rho_{P_{\beta}}$. 
\end{itemize}
\end{preu}

As a consequence of Lemma \ref{Plem}, the radius $ \rho_{P_{\beta}^*} = \frac{1}{|\beta|}$ is the unique solution of the equation $\sum\limits_{x\in P_{\beta}}z^{l(x)} = 1$.

\begin{defi}
 Let $ X $ be the prefix code of a symbolic dynamical system $ S $ with language $ L_S $. A message of $ L_S$ is a word $ A_1 A_2 \cdots A_k $ where each word $ A_i$ belongs to the code $ X $.
\end{defi}

\begin{prop}\label{MixProp 1}
 Let $ \beta < -1 $ and $ (d_i)_{i \geq 1} $ be the $\beta$-expansion of $ l_{\beta} = \frac{\beta}{1-\beta}$. We assume that $ B_n $ counts the number of messages of length $n$ of the code of the maximal measure support. Then, there is a constant $ C $ such that $B_n< C |\beta|^n$.
\end{prop}

\begin{preu}
 %\begin{itemize}
 % \item[(i)]
  Suppose $\beta < - 1 $ and let $ H_n$ counts words of length $ n $ in $ L_{\beta}$. We have $ B_n < H_n$. But:
  \begin{equation*}
  H_n = \sum\limits_{k=1}^n(-1)^k (d_{k-1}-d_k)H_{n-k}+1.
  \end{equation*}
And thus,
\begin{equation*}
 \sum\limits_{n \geq 0} H_nz^n = \dfrac{\sum\limits_{n \geq 0}z^n}{1-\sum\limits_{n \geq 1}(-1)^n(d_{n-1}-d_n)z^n}.
\end{equation*}
The sequence $(d_n)_{n \geq 1}$ is bounded, $ 1-\sum\limits_{n \geq 1}(-1)^n(d_{n-1}-d_n)z^n $ is then analytic. Thus, its roots are isolated points. We know that $ \frac{1}{|\beta|}$ is the smallest pole in modulus of $ H= \sum\limits_{n \geq 0} H_nz^n$. This is a single root of $ 1-\sum\limits_{n \geq 1}(-1)^n(d_{n-1}-d_n)z^n $ and no other root in the disk of center 0 and radius $ \frac{1}{|\beta|}$. Thus, in this disk, we can write $ H $ as follows:
\begin{equation*}
 H = \frac{\alpha}{1-|\beta|z} + \frac{P}{Q}
\end{equation*}
where, in the disk, $ P$ is analytic and $ Q $ has no root. Let $ \sum\limits_{n \geq 0}a_n z^n $ be the power series expansion of $\frac{P}{Q}$. Let $ R = \sup\{r , (a_nr^n)_{n \geq 0} \text{ bounded }\}$ be its radius. We have $ \frac{1}{|\beta|} \leq R < +\infty$. There is $ K > 0$ such that $R^n|a_n| < K$. That means that $ |a_n | < \frac{K}{R^n}< K|\beta|^n$. By using the formal power expansion of $ \frac{\alpha}{1-|\beta|z} $, we deduce that
\begin{equation*}
 H_n = \alpha |\beta|^n + a_n \leq |\alpha | |\beta|^n + K|\beta|^n < C |\beta|^n.
\end{equation*}
\end{preu}

Let $ \beta $ be a real number smaller than $-1$. We recall that $ I_{\beta} = [\frac{\beta}{1-\beta}, \frac{1}{1-\beta} )$ and the map $ T_{\beta} $ as defined in equation (\ref{T}).  In \cite{MR2974214}, the authors proved the intrinsic ergodicity of the dynamical system $ (I_{\beta}, T_{\beta} )$. 
This question has not been addressed for  the negative $\beta$-shift $ S_{\beta}$. But, due to the conjugacy of $ (I_{\beta}, T_{\beta}) $ and $ (S_{\beta}^r, \sigma) $, the right one-sided $\beta$-shift is an intrinsically ergodic system. The aim of the following subsection is to prove the intrinsic ergodicity of the two-sided negative $ \beta$-shift by using techniques of \cite{MR858689}. To do that, we need the prefix recurrent code $ P_{\beta}$ to construct an associated tower.  If $ \beta< -\gamma_0$, $ P_{\beta} = \C_{\beta} $, code of the sequences which end on the form $ d_1 \cdots d_k j $ with $ (-1)^{k+1}(d_k-j)< 0 $ and $ j \neq d_1$. If $ -\gamma_n \leq \beta < -\gamma_{n+1}$, $ P_{\beta} = \Delta_{n, \beta}$.

\subsection{Intrinsic ergodicity of the two-sided negative $ \beta$-shift}

The existence of the positive recurrent prefix code $P_{\beta}$ implies the uniqueness of the measure $\overline{\mu}$ with entropy $\log |\beta|$ on the tower $(\Omega, T)$ attached to $P_{\beta}$. This measure induces the Bernoulli probability $\mu$ on the free monoid generated by $P_{\beta}$ and defined by
\begin{equation}
 \mu([x]) = \frac{1}{|\beta|^{l(x)}}, \text{ $ x\in P_{\beta} $ }.
\end{equation}
Let $\P(S_\beta) = \{ A | A\subset S_{\beta} \}$ and $\nu$ be the map from $\P(S_\beta)$ into $[0, 1]$ which coincides with $\mu$ on $P_{\beta}^{\Z}$ and zero on all subset of the complement of $P_{\beta}^{\Z}$.
\begin{equation}
 \nu (B) = \begin{cases}
        \mu (B) & \text{ if $ B \subset P_{\beta}^{\Z} $ }  \\
        0       & \text{ if $ B \subset S_{\beta}\setminus P_{\beta}^{\Z}$ }
       \end{cases}
\end{equation}
%The set $S_{\beta}\setminus W(P_{\beta})$ consists of sequences which end by $d(l_{\beta}, \beta)$ and the free monoid generated by $D$.
We are going to prove that a measure of maximal entropy cannot be carried by a subset of $ S_{\beta}\setminus P_{\beta}^{\Z}$.

Let $ (\Omega, T) $ be the tower of $ P_{\beta} $ and denote by $ f $ the map from $ (\Omega, T) $ into $\A^{\Z}$ such that:
\begin{equation}
 f( (x_n)_{n \in \Z}, i) =  (y_n)_{n \in \Z},
 \label{(2)}
\end{equation}
where $ x_k \in P_{\beta} $, $ (y_i)_{i \in \Z} $ is the unique sequence of $ \A= \{0, 1, \cdots, d_1 \} $ such that for any integer $ n $,
\begin{equation*}
x_{-n}x_{-n+1}\cdots x_{-1}x_0x_1\cdots x_n = y_{-p}y_{-p+1}\cdots y_{-1}y_0y_1\cdots y_{q-1}
\end{equation*}
with
\begin{equation*}
p+q = l(x_{-n}x_{-n+1}\cdots x_{-1}x_0x_1\cdots x_n)
\end{equation*}
and
\begin{equation*}
p=l(x_{-n}x_{-n+1}\cdots x_{-1})+i-1.
\end{equation*}
That is $ y_0 $ is the $i$-th letter of $x_0$. %Indeed, let $ x $ be the two-sided admissible sequence $x = \cdots z_{-k}{z_{-k+1}\cdots z_1z_0\bullet z_1\cdots z_k z_{k+1}\cdots $

\begin{exple}
 Let $\beta$ such that $d(l_{\beta}, \beta) = \bullet2 (1)^{\infty}$. We set $D=\{2(11)^k, k \in \N\}$
 \begin{equation}
  \C_{\beta} = \{0, 1,  x10 | x \in D \}.
 \end{equation}
We set $x_{n}=210$ for $n\leq -1$, $x_0=21110$ and $x_n = 2210$. Then

\begin{equation*}
 f((x_n)_{n\in \Z}, 2) = \cdots 210 \hspace{1mm} 210\hspace{1mm} 210 \hspace{1mm} 21 \bullet 110 \hspace{1mm}\overline{2210}.
\end{equation*}

\end{exple}

\begin{rem}
 The map $ f $ is defined in \cite{MR858689}. As shown in \cite{MR858689}, $ f\circ T = \sigma \circ f $ and $f$ is one to one. Let $(z_n)_{n \in \Z}\in \A^{\Z}$, $(x_n)_{n\in \Z}\in P_{\beta}^{\Z}$ such that:
 \begin{equation*}
\cdots x_{-m}\cdots x_{-1}x_0 x_1 x_2 \cdots x_m \cdots = \cdots z_{-n}\cdots z_{-1} z_0 \bullet z_1 \cdots z_n \cdots
\end{equation*}
where $z_0$ is the $i$-th letter of $x_0$. The inverse map $f^{-1}$ of f from $P_{\beta}^{\Z}$ to $\Omega$ is given by

\begin{equation}
 f^{-1}((z_n)_{n \in \Z}) = ((x_n)_{n \in \Z}, i).
\end{equation}
\end{rem}
 Each $\sigma$-invariant measure $\pi$ on $P_{\beta}^{\Z}$ generates a measure $\pi \circ f$ on $\Omega$ with the same entropy (see Proposition 2.17 of \cite{MR858689}). Thus the uniqueness of the measure of the entropy $\log |\beta|$ on $\Omega$ implies that there is a unique measure of same entropy on $P_{\beta}^{\Z}$. 

\subsubsection{Hybrid admissible sequences}

\begin{defi}
In what follows, an admissible sequence $ (x_i)_{i \in \Z}$ shall be said to be hybrid if there is an integer $ k \in \Z $ such that one of the following conditions hold:
\begin{itemize}
\item[(a)] $ (x_i)_{i \geq k } \in \C_{\beta}^{\N} \text{ and } (x_i)_{ i < k } \not\in \C_{\beta}^{\N} $;
\item[(b)] $ (x_i)_{i \geq k } \not\in \C_{\beta}^{\N} \text{ and } (x_i)_{ i < k } \in \C_{\beta}^{\N} $;
\item[(c)] $(x_i)_{i \geq k } \in \Delta_{j,\beta}^{\N} \text{ and } (x_i)_{ i < k } \not\in \Delta_{j, \beta}^{\N} $;
\item[(d)] $(x_i)_{i \geq k } \not\in \Delta_{j,\beta}^{\N} \text{ and } (x_i)_{ i < k } \in \Delta_{j, \beta}^{\N}$.
\end{itemize}
\end{defi}

\begin{exple}\label{Ex5}
We suppose that $ d_1d_2 \cdots  = 2 (1)^{\infty} $, $210 \in \C_{\beta}$. The following sequences are hybrid: $ \cdots 210 210 210 \cdots 210 \bullet (2)^{\infty}$, $ \cdots 2 222 \cdots 2 \bullet (210)^{\infty}$, $ \cdots 210 210 \cdots 210 \bullet 2 (1)^{\infty} $, $\cdots 211 211 \cdots 211 \bullet 2(1)^{\infty} $. 
\end{exple}

Let $ N_{\beta, k} $ be the subset of $ N_{\beta} $ of sequences $ (x_i)_{i \geq \Z} $ such that one the following situations occurs: 
\begin{itemize}
\item[(e)] $ (x_i)_{i \leq k} \not\in \C_{\beta}^{\N} $, there is $  m $ such that $ x_{k+1} \cdots x_{k+m-1} = d_1\cdots d_{m-1} $ and $(-1)^m(d_m - x_{k+m}) < 0 $, $ x_{k+m}\neq d_1$ and $ (x_i)_{i \geq k+1} \in \C_{\beta}^{\N}$, 
\item[(f)] $ (x_i)_{i \leq k} \in \C_{\beta}^{\N} $ and $ (x_i)_{i \geq k+1} \not\in \C_{\beta}^{\N}$; 
\item[(g)]  $ (x_i)_{i \leq k} \not\in \Delta_{i,\beta}^{\N} $, there is $  m $ such that $ x_{k+1} \cdots x_{k+m-1} $ is an initial segment which belongs to  $\Delta_{i,\beta}$ and $ (x_i)_{i \geq k+1} \in \Delta_{i,\beta}^{\N}$, 
\item[(h)] $ (x_i)_{i \leq k} \in \Delta_{i, \beta}^{\N} $ and $ (x_i)_{i \geq k+1} \not\in \Delta_{i, \beta}^{\N}$. 
\end{itemize}
The sequences given in Example \ref{Ex5} belong to $ N_{\beta, 0}$. 
\begin{prop}\label{N}
The set $N_{\beta} $ of hybrid sequences  is a null set with respect to any $ \sigma$-invariant measure on $ S_{\beta}$.
%Denote by $ N $, the subset of admissible sequences $ (x_i)_{i \in \Z}$ of the two-sided negative $ \beta$-shift for which there is an integer $ k \in \Z $ such that the sub-sequence: 
%\begin{equation}
%\begin{cases}
%(x_i)_{i \geq k } \in \C_{\beta}^{\N} &\text{ and } (x_i)_{ i < k } \not\in \C_{\beta}^{\N} \text{ or } \\
%(x_i)_{i \geq k } \in \Delta_{j,\beta}^{\N} &\text{ and } (x_i)_{ i < k } \not\in \Delta_{j, \beta}^{\N}
%\end{cases}\label{71}
%\end{equation}
%for some $ j $, is negligible with respect to any $ \sigma$-invariant measure on $ S_{\beta}$.
\end{prop}

\begin{preu}

We have $ N_{\beta} = \underset{k \in \Z}{\bigcup} N_{\beta,k} $. It is easy to show that $ \sigma^{-1} N_{\beta,k} = N_{\beta, k-1} $ and $ N_{\beta,k} \cap N_{\beta, k-1} = \emptyset$. Thus, for a $ \sigma$-invariant measure $ \mu $, we have $ \mu (N_{\beta, k}) = \mu (N_{\beta, k-1}) $ and $ \mu (N_{\beta}) = \sum\limits_{k \in \Z} \mu (N_{\beta, k})$. If $ \mu(N_{\beta, k}) = \mu(N_{\beta,0}) \neq 0$, then $ \mu(N_{\beta})=+\infty $. Hence the result follows. 
\end{preu}

An admissible sequence $(x_i)_{i\in \Z}$ which is not in $ N_{\beta}$ belongs to $ \C_{\beta}^{\Z}$, $ \Delta_{\beta, j}^{\Z} $ for some integer $ j \in \N$, or $ \exists m \in \Z$: $ \forall j < m$, $ \exists k \leq j$ with $ x_{k+1}\cdots x_m = d_1 \cdots d_{m-k} $. That is $ x \in \underset{ j \in \Z}{\bigcup} \sigma^{j} \underset{ l \geq 1}{\bigcap}\underset{i \geq l}{\bigcup}\sigma^{-i} [d_1 \cdots d_i ] $. Let $ \Lambda_{\beta} $ be such of admissible words. %the subset of admissible sequences $ (x_i)_{i \in \Z}$ such that $ \forall j < m$, $ \exists k \leq j$ with $ x_{k+1}\cdots x_m = d_1 \cdots d_{m-k} $. 

\begin{equation}
\Lambda_{\beta} = \underset{ j \in \Z}{\bigcup} \sigma^{j} \underset{ l \geq 1}{\bigcap}\underset{i \geq l}{\bigcup}\sigma^{-i} [d_1 \cdots d_i ]. 
\end{equation}

For this last class of admissible sequences, it is possible to push infinitely to the left the characteristic sequence $ (d_i)_{i \geq 1}$. This case has been highlighted in \cite{hofbauer1978beta} for $ \beta > 1 $. For such a $ \beta$, this case occurs when $ \lim\limits_{n \rightarrow +\infty} \sigma^n(d_{\beta}(1)) = d_{\beta}(1) $ (in our case, $ \beta < -1$, similarly, $ \lim\limits_{n \rightarrow +\infty} \sigma^n(d(l_{\beta}, \beta)) = d(l_{\beta}, \beta) $).

\begin{equation}
S_{\beta} = N_{\beta} \cup \Lambda_{\beta} \cup \C_{\beta}^{\Z} \cup  \left(  \underset{ j \in \N}{\bigcup} \Delta_{j,\beta}^{\Z}\right)\label{eqAB} .
\end{equation}

\subsubsection{If $ \beta< - \gamma_0 $}

\begin{rem}
For a given initial segment with odd length $ d_1 \cdots d_{2n+1} $, the infinite sequence $ \overline{d_1 \cdots d_{2n+1}} $ is not admissible if and only if $ 2n_i \leq 2n+1 < 2n_i+p_i-1 $.  
\end{rem}

\begin{prop}\label{Jur}
Let $(d_i)_{i \geq 1}$ be the characteristic sequence defined in \ref{(D)} and \ref{(6)}. For all $ k \in \N^{*}$ there is $ p \in \N^{*} $ such that $ 2p+1 \geq k $ and $ \overline{d_1\cdots d_{2p+1}} $ is admissible.
\end{prop}

\begin{preu}
We consider the characteristic sequence on the general form
\begin{equation*} 
d_1 \cdots d_{2n_1-1} d_1 \cdots d_{p_1} d_{2n_1+p_1} \cdots d_{2n_2-1}d_1\cdots d_{p_2} d_{2n_2+p_2}\cdots d_{2n_3-1} \cdots.
\end{equation*} 
For all integer $ i \in \N^{*} $, $ \overline{d_1 \cdots d_{2n_i-1}}$ is admissible. Thus, if there is an infinite number of $ i $ such that $ d_1 \cdots d_{2n_i-1}$ is followed by $ d_1 \cdots d_{p_i}$ in the characteristic sequence, for all $ k \in \N^{*}$ there is $ i \in \N^{*} $ such that $ 2n_i-1 \geq k $ and $ \overline{d_1\cdots d_{2n_i-1}} $ is admissible. Hence the result follows. 

If there is a finite number of $ i $ such that $ d_1 \cdots d_{2n_i-1}$ is followed by $ d_1 \cdots d_{p_i}$, we set $i_0$ the maximum of such of integers. Since $ p_i <2n_i-1$ because the characteristic sequence is not periodic with odd period, for all $ k $ such that $ 2k+1 \geq 2n_i+p_i $, $ \overline{d_1\cdots d_{2k+1}} $ is admissible. Hence the result follows. 
\end{preu}

As a consequence of Proposition \ref{Jur}, all initial segment is a word of the language of $ \Delta^{\Z}\cap S_{\beta}$, where $ \Delta $ is the set of initial segment of odd length. 

\begin{prop}\label{Prop G}
Let $ \beta < - \gamma_0 $. We assume that $ \mu $ is an ergodic measure of maximal entropy on $ S_{\beta}$. Then $ \mu(\Lambda_{\beta}) = 0$.
\end{prop}
The proof of Proposition \ref{Prop G} is based on Theorem 1 of \cite{hofbauer1978beta}.

\begin{preu}
Let $ \Lambda_{\beta,m} $ be the set of sequences $ x=(x_i)_{i \in \N}$ such that: $ \forall j < m $, $ \exists k \leq j$, $ x_k \cdots x_{m-1} = d_1 \cdots d_{m-k} $. 
\begin{equation}
\Lambda_{\beta, m} = \sigma^m \underset{l \geq 1}{\bigcap} \left(  \underset{i \geq l }{\bigcup} \sigma^{-i}[d_1 \cdots d_{i-1} ]\right).
\end{equation}
We have $\Lambda_{\beta} = \underset{ m \in \Z}{\bigcup} \Lambda_{\beta,m} $, $ \Lambda_{\beta, m} = \sigma^{-1} \Lambda_{\beta, m+1} $. For an invariant measure $ \mu $, $ \mu (\Lambda_{\beta, m} )= \mu(\Lambda_{\beta, m+1}) $. And thus, $ \mu (\Lambda_{\beta} ) = \mu \left(  \underset{m \in \Z}{\bigcap} \Lambda_{\beta,m} \right) $.
\begin{equation}
\underset{m \in \Z}{\bigcap} \Lambda_{\beta,m} = \{ x \in S_{\beta}: \forall k  \in \N, \exists r, s, \hspace{1mm} r<-k, s > k \text{ such that } x_rx_{r+1} \cdots x_s = d_1 \cdots d_{s-r+1} \}.
\end{equation}
According to the equality above, observe that for all $ x = (x_i)_{i \in \Z}\in \underset{m \in \Z}{\bigcap} \Lambda_{\beta,m}$, for all integers $ r< s $, $ x_r \cdots x_s $ is a string which appear in an initial segment. That is, the language of $\underset{m \in \Z}{\bigcap} \Lambda_{\beta,m}$ is included in the language of $ \Delta^{\Z}\cap S_{\beta} $. But, the entropy of $ \Delta^{\Z}\cap S_{\beta} $ is smaller than $ \log \gamma_0 $. Since $ \Lambda_{\beta} $ is a $\sigma$-invariant subset of $S_{\beta}$, it follows that it is a null set with respect to all measure of maximal entropy. 
\end{preu}

As a consequence, we have the following corollary:

\begin{coro}\label{Support}
Let $ \beta $ be a real number smaller than $ -\gamma_0$. Then a measure of maximal entropy on $ S_{\beta}$ is carried by $ \C^{\Z}$.
\end{coro}

\begin{preu}
The set $ \underset{ j \in \N}{\bigcup} \Delta_{j, \beta}^{\Z} $ is $ \sigma $-invariant. Endowed with the shift, it has a topological entropy less $ \log \gamma_0 $ since it is included in the set of admissible concatenations of initial segments of odd length. Thus, it can not be the support of a measure of maximal entropy. As a consequence, a measure of maximal entropy on $ S_{\beta}$ is carried by $ \C^{\Z} $. 
\end{preu}
%As consequence of the previous proposition, 

\subsubsection{If $ - \gamma_0 < \beta < - 1 $ }

Suppose $ - \gamma_0 < \beta < -1 $. There is $ n \in \N $ such that  $ -\gamma_n < \beta \leq -\gamma_{n+1} $. Moreover, we can find a unique real number $ x \in ] -\infty, - \gamma_0 ]$ such that $ d(l_{\beta}, \beta) = \phi^{n+1} (d(l_x, x)) $. 
So, $ \phi^{n+1}(S_x) $ is a sub-set of $ S_{\beta}$ coded by $ \phi^{n+1} (\C_{x})$. A sequence $ (x_i)_{i \in \Z}\not\in \phi^{n+1} S_x $ belongs to $ \Lambda_{\beta}$, $ N_{\beta}$, $\C_{\beta}^{\Z} $ or $ \Delta_{j, \beta}^{\Z}$ for some integer $j $. If $ j\geq n $, it is easy to see that $ \Delta_{j, \beta}^{\Z} \subset \phi^{n+1} S_x$. Furthermore, in the previous section, we have shown that for any $\sigma$-invariant measure $ \mu$ on $ S_{\beta}$, $ \mu(\Lambda_{\beta})= \mu \left(  \underset{m \in \Z}{\bigcap} \Lambda_{\beta,m} \right) $ and thus $\mu \left( \Lambda_{\beta}\setminus \underset{m \in \Z}{\bigcap} \Lambda_{\beta,m} \right) = 0 $ . Moreover, we have seen that for all sequence $ (x_i)_{i \in \Z}$ taken in $\underset{m \in \Z}{\bigcap} \Lambda_{\beta,m}$, for $ k < n $, $ x_k x_{k+1} \cdots x_n $ appears in an initial segment. That is, $ x_k \cdots x_n $ appears in $ d(l_{\beta}, \beta)= \phi^{n+1}(d(l_x, x))$. This implies that $\underset{m \in \Z}{\bigcap} \Lambda_{\beta,m} \subset \phi^{n+1}(S_x)$. Thus, \ref{eqAB} can be written as follows: 

\begin{equation}
S_{\beta} = \phi^{n+1}(S_x) \cup N_{\beta} \cup \C_{\beta}^{\Z} \cup  \left( \Lambda_{\beta}\setminus \underset{m \in \Z}{\bigcap} \Lambda_{\beta,m}\right)  \cup \left( \underset{k=0}{\overset{n-1}{\bigcup}}\Delta_{\beta, k}^{\Z}\right).
\end{equation}

%Remark that:

%\begin{equation}
%S_{\beta} = \phi^{n+1}(S_x) \cup N \cup \left( \underset{k=0}{\overset{n-1}{\bigcup}}\Delta_{\beta, k}^{\Z}\right)
%\end{equation}
%where $ N $ is described in the proof of Proposition \ref{N}, $ N = \underset{k \in \Z}{\bigcup} N_k $, with 
%\begin{equation}
%  N_k = \lbrace (x_i)_{i \in \Z} \in S_{\beta} \mid k = \max\{m, (x_i)_{i< m} \not\in \Delta_{j, \beta}^{\N} \text{ and } (x_i)_{i \geq m} \in \Delta_{j, \beta}^{\N}\} \rbrace.
%\end{equation}
\begin{coro}\label{support 2}
Let $ \beta$ be a real number such that $ -\gamma_n \leq \beta < -\gamma_{n+1}$. A measure of maximal entropy is supported by $\Delta_{n, \beta}^{\Z} = \phi^{n+1} (\C_{x})^{\Z}$, for some $ x \in ]-\infty, -\gamma_0 [$. 
\end{coro}

\begin{preu}
The set of admissible hybrid sequences is a null set with respect to any $ \sigma$-invariant measure. If $ \beta \in ] -\gamma_n, -\gamma_{n+1} [ $, then $ \C_{\beta} = \{0\}$ and $ \Delta_{\beta, k} = \{u_k\} $. This implies that the topological entropy of $\Delta_{\beta, k}^{\Z} $, just like that of $ \C_{\beta}^{\Z}$, endowed with the shift is 0 for any nonnegative integer $ k $ smaller than $ n-1$.  Thus, the support of a measure of maximal entropy is included in $ \phi^{n+1}(S_x) $. 

It is easy to see that the map $ \phi $ is an injective morphism. Thus, $ \phi^{n+1} $ establishes a bijection from $ S_x $ to $ \phi^{n+1} (S_{x})$.

\begin{equation*}
\phi^{n+1} (S_{x}) = \phi^{n+1}(N_x )\cup \phi^{n+1}(\Lambda_x) \cup \phi^{n+1}(\C_x)^{\Z} \cup  \left(  \underset{ j \in \N}{\bigcup} \phi^{n+1}(\Delta_{j,x})^{\Z}\right) .
\end{equation*}
 Furthermore, $ \phi^{n+1}$ operates a change of alphabet from $\{0, 1, \cdots, \rfloor \vert x \vert\lfloor \} $ to $ \{u_n, u_nu_{n-1}u_{n-1}, \cdots, u_n(u_{n-1}u_{n-1})^{\rfloor \vert x \vert\lfloor}$. Moreover, $ \phi^{n+1}(S_x) $ is coded by $ \phi^{n+1} (\C_{\beta})^{\Z}$.  
From Corollary \ref{Support}, a measure of maximal entropy on $ \phi^{n+1}(S_x) $ is carried by $ \phi^{n+1} (\C_{\beta})^{\Z}$.

If there are two measures of maximal entropy $\log |\beta|$, there are carried by the same support coded by $P_{\beta}$. Then, they coincide on subsets of $P_{\beta}^{\Z}$ and subsets of $S_{\beta}\setminus P_{\beta}^{\Z}$.That is the measure of maximal entropy is unique.
\end{preu}

\section{Ergodic measure and Proofs of Theorems \ref{Th tau}, \ref{Th surp}, \ref{th3} and \ref{Theo 4}}

From now on, we denote by $ P_{\beta} $ the code of the support of the measure of maximal entropy.
\begin{equation}
P_{\beta} = \begin{cases}
    \C_{\beta} &\text{ if $ \beta < - \gamma_0 $ } \\
		\Delta_{n,\beta} &\text{ if $ -\gamma_{n}\leq \beta < -\gamma_{n+1}$}.
		\end{cases}
\end{equation}

\begin{lem}\label{last}

 Let $x$ and $ \beta $ be two real numbers such that $ x \in ] -\infty, - \gamma_0 ]$, $ \beta \in ] -\gamma_n, -\gamma_{n+1} ]$ and $d^{*}(l_{\beta}, \beta) = \phi^{n+1} (d^{*}(l_x, x))$. Then, in the meaning of the alternating order
 \begin{equation}
 u_nd^{*}(l_{\beta}, \beta) = \sup \phi^{n}(\Delta_{n, \beta}^{\N})=\sup \phi^{n+1}(S_x^d) \label{sup}
  \end{equation}
and it is a $\beta$-expansion of

 \begin{equation*}
  t_n = \prod\limits_{k = -1}^{n-1}(1+\frac{1}{\beta^{l(u_k)}}) - \frac{\beta^{l(u_n)}-2}{\beta^{l(u_n)-1}(\beta-1)}.
 \end{equation*}

\end{lem}

In what follows, given a sequence $ (a_n)_{n \in \N}$ (finite or infinite), $(a_1a_2\cdots)(-z) = \sum\limits_{n \geq 1} a_n z^n$. So, if $ (a_n)_{n \in \N}$ belongs to $S_{\beta}^d$, then this sequence is a $\beta$-representation of $ (a_1a_2 \cdots)(-\frac{1}{\beta})$.

\begin{preu}\textbf{of Lemma \ref{last}}
Note that $\sup \phi^{n}(\Delta_{n, \beta}^{\N})=u_nd^{*}(l_{\beta}, \beta) = \sup \phi^{n+1}(S_x^d)$ thanks to the fact that, in the meaning of the alternating order, $ \phi$ is an increasing map.
Denote by $ b_1b_2 \cdots b_{l(u_n)} $, the word on $ \{1, 2\}^*$ such that if $ u_n = \phi^n(1) = \epsilon_1\epsilon_2 \cdots \epsilon_{l(u_n)}$, $ b_i = \epsilon_i + 1$. From \cite{NguemaNdong20161},

\begin{equation}
 (1b_1b_2\cdots b_{l(u_n)})(-z) = \prod\limits_{k=-1}^{n-1}(1+z^{l(u_k)}) + z^{l(u_n)}.
\end{equation}
And
\begin{equation}
(1b_1b_2\cdots b_{l(u_n)})(-z) = (0u_n)(-z)+\frac{1-z^{l(u_n)}}{1-z}.
\end{equation}
Thus
\begin{equation}
\begin{aligned}
 (0u_nd(l_{\beta}, \beta))(-\frac{1}{\beta}) &= (0u_n)(-\frac{1}{\beta}) + \frac{1}{\beta^{l(u_n)}}l_{\beta} \\
                                             &= \prod\limits_{k = -1}^{n-1}(1+\frac{1}{\beta^{l(u_k)}}) - \frac{\beta^{l(u_n)}-2}{\beta^{l(u_n)-1}(\beta-1)}.
 \end{aligned}
\end{equation}
\end{preu}

The equation \eqref{sup} implies that $ \phi^{n+1}(S_x)$ is the closure of $ \Delta_{n, \beta}^{\Z}$.

\vspace{0.5cm}
\begin{preu}\textbf{of Theorem \ref{Th surp}}
The map $ \Upsilon$ is well defined. Indeed, for any $ x \in ] -\gamma_0, - 1 [$, since $(\gamma_n)_{n \geq 0}$ decreases and its limit is 1, there is a unique integer $ n \in \N $ such that $ x\in ]-\gamma_n, -\gamma_{n+1} ]$. From Corollary \ref{copp2}, there exists a unique real number $\beta \in ]-\infty, -\gamma_0 ]$ such that 
\begin{equation*}
d^{*}(l_x, x) = \phi^{n+1}(d^{*}(l_{\beta}, \beta)).
\end{equation*}
That is, $ \beta = \Upsilon(x)$. Thus, to each number of the interval $ ]-\gamma_0, -1[$, we associate a unique element $ \beta$ of $]-\infty, -\gamma_0]$.

 The proof of Theorem \ref{Th surp} is obtained by using Proposition \ref{pp2}, Corollary \ref{copp2} and Remark \ref{remg}.
For all $x\in ]-\infty, -\gamma_0 ]$, there is a sequence $(\beta_n)_{n \in \N} $ such that $\beta_n \in ] -\gamma_n, -\gamma_{n+1}]$ and $d^{*}(l_{\beta_n}, \beta_n) = \phi^{n+1} (d^{*}(l_x, x))$, that is $x=\Upsilon(\beta_n)$ for any $n\geq 0$.
\end{preu}

It is well known that the map $x\mapsto d(l_x, x)$ is increasing. The same holds for the map $x\mapsto d^{*}(l_x, x)$. This implies that $\Upsilon$ is an increasing map. The morphism $\phi$ changes an alphabet $\{0, \cdots, k\}$ to $\{1, 100, 10000, \cdots, 1(00)^k\}$, with $ k = \rfloor \vert x\vert \lfloor$. However $ \phi $ does not change the cardinality of the alphabet. The string $(00)^k$ is the longest sequence of zeros to add after $1$ in an admissible string. Note that $d^{*}(l_{\beta}, \beta) $ is a sequence of $ \{0, 1, \cdots, \rfloor |\beta| \lfloor\}$, where $\rfloor |\beta| \lfloor$ denotes the largest integer strictly smaller than $|\beta|$. If $x$ tends to $-\gamma_0$, the size of the longest string of zeros after 1 tends to $+\infty$. Thus $\rfloor |\Upsilon(x)| \lfloor$ tends to $ + \infty$. And then $\lim\limits_{x\rightarrow -\gamma_0} \Upsilon(x) = - \infty$.

\begin{preu}\textbf{of Theorem \ref{Th tau}}
 Let $\beta\in ] -\gamma_0, -1 [$, from Theorem \ref{Th surp}, there is $x\in ]-\infty, -\gamma_0 ]$ such that $ d^{*}(l_{\beta}, \beta) = \phi^{n+1} (d^{*}(l_{x}, x))$ if $\beta\in ]-\gamma_n, - \gamma_{n+1} ]$. We set $g = \phi^{n+1}$. From Proposition 2.10 of \cite{MR858689} and Remark \ref{remg}, $g(S_x)$ is coded by $g(\C_x) = \Delta_{n, \beta}$. The result follows by the fact that if $\beta\in]-\gamma_n, -\gamma_{n+1}[$, $\Delta_{n, \beta}$ codes the support of the measure of maximal entropy and by using Lemma \ref{last}.

 %To prove that the map $ \phi^{n+1} $ transforms the $ x $-expansions of real numbers of the interval $ [l_x, r_x [$ to $\beta$-expansions of real numbers of $ [l_{\beta}, t_n [$, remark at first that in the meaning of the alternating order, $ \sup \Delta_{n, \beta}^{\N} = u_nd(l_{\beta}, \beta) $.
\end{preu}

\begin{theo}\label{Mix}
 The measure of maximal entropy on the negative beta-shift is mixing.
 \end{theo}
Before proving Theorem \ref{Mix}, let us show the following result:

\begin{prop}\label{GCD}
The $\gcd$ of lengths of words of codes previously constructed is 1.
\end{prop}

\begin{preu}
 For $ \beta \leq -\gamma_0$, the $\beta$-shift is coded (by $ \C_{\beta} $ if the inequality is strict and by $ \{1, 00\} $ if $ \beta = \gamma_0$). And also, the code contains at least one word of length 1.

 If $ \beta \in [ -\gamma_0, -\gamma_1 )$, the support is coded by $\Delta_{0, \beta}$ which contains $ 1 = d_1 $.

 Therefore, consider $ \beta $ such that $ -\gamma_{n+1} > \beta \geq -\gamma_n $ with $ n > 1$. In this case, the support of maximal entropy measure is coded by $ \Delta_{n, \beta}$. The words of this set are of the form
 \begin{equation}
 u_nv_n^{n_1+1}u_{n}v_n^{n_2}\cdots u_n v_n^{n_{2k}}u_nv_n^t,
 \end{equation}
 with $ 0 \leq n_{2k+1}-1$ and $ 0 \leq k$.

If $ n_1 \neq 0$, $ u_n $ and $ u_nv_n $ belong to $ \Delta_{n, \beta} $. The integer $ l( u_n ) $ and $ l( v_n ) $ are relatively prime since
 \begin{equation}
 l( u_n ) = l( v_n ) - (-1)^n.
 \end{equation}
Thus, $ l( u_n ) $ and $ l( u_n v_n ) $ are relatively prime too.

Note that $ v_n^{n_1+1}$ is the longest sequence of $ v_n $ in the support of a measure of maximal entropy. Thus, if $ n_1 = 0$, $n_2=0$, $ u_n$ and $ u_n v_n u_nu_n$ belong to the code. But $ l( u_n ) $ and $ l( u_n v_n u_nu_n )$ are relatively prime. It follows that, for all $ \beta <- 1$, the $\gcd$ of lengths of words belonging to the code of the support of maximal entropy measure is 1. %That is the measure of maximal entropy is mixing.
\end{preu}

\begin{preu}\textbf{ Theorem \ref{Mix}}\\
An immediate consequence of the Proposition \ref{GCD} is that the restriction of the measure of maximal entropy on its support is mixing. Note that, if $ x $ is an intransitive word, $ [x] $ is $\sigma$-invariant. And then, for all $ n $ and $ y $ in the code of support,
\begin{equation}
\sigma^{-n}[x]\cap [y] = \varnothing
\end{equation}
Thus
\begin{equation}
\lim\limits_{n\rightarrow +\infty} \mu(\sigma^{-n}[x]\cap [y]) = 0 = \mu([x])\mu([y])\label{eqmix}
\end{equation}
since $\mu([x]) = 0 $. Moreover, for all $ n $, $\sigma^{-n}[y]\cap [x] \subset [x ]$ and then
\begin{equation}
\lim\limits_{n\rightarrow +\infty} \mu(\sigma^{-n}[y]\cap [x]) = 0 = \mu([y])\mu([x]).
\end{equation}
If now, $ x $ and $ y $ are both intransitive words, $\sigma^{-n}[x] \cap [y]$ is a null set with respect to the measure of maximal entropy. Then \eqref{eqmix} is also satisfied. This proves Theorem \ref{Mix}.
\end{preu}

\begin{rem}\label{remMix}
 We have proved that
 \begin{equation*}
  \begin{cases}
   \rho_{P^*_{\beta}} < \rho_{P_{\beta}} \\
   P_{\beta} \text{ prefix } \\
   \gcd(\{l(x) | x \in P_{\beta}\}) =  1.
  \end{cases}
 \end{equation*}
Denote by $ e_{k,i}$ the $i$-th words of length $k$ in $P_{\beta}$. From \cite{MR939059}, the sequence $ (e_{k,i})_{k,i \geq 1}$ is generic for a measure $\mu_{P_{\beta}}$ (Champernowne measure) which has entropy $-\log \rho^*_{P_{\beta}}=\log |\beta|$. This measure is strongly mixing. Since $ S_{\beta} $ is intrinsically ergodic, it follows that, $\mu_{P_{\beta}}$ is the unique measure of maximal entropy.
\end{rem}

\begin{preu} \textbf{ of Theorem \ref{th3}} \\
If $ \beta < -\gamma_0$, $\C_{\beta} $ codes $ S_{\beta}$ (see \cite{FloarXiv,NguemaNdong2019}). From Theorem \ref{Th tau}, Lemma \ref{lm2}, 6 and 7 of \cite{NguemaNdong2019}, the unique measure of maximal entropy is supported by a set coded by $P_{\beta}$ which is positive recurrent from Lemma \ref{Plem}.

Any invariant probability $ \nu $ on $(P_{\beta}^{\Z}, \sigma_{P_{\beta}})$ with finite average length $L(P_{\beta}, \nu) $ is induced by a unique invariant probability measure $ \overline{\mu} $ of $(\Omega, T ) $ (see \cite{MR858689}). The link between the entropies of the two measures is given by the Abramov formula:
\begin{equation}
h(\overline{\mu})L(\nu, P_{\beta}) = h(\nu).\label{abra}
\end{equation} 
and $\overline{\mu}\circ f^{-1}$ is a probability on $(A^{\Z}, \sigma)$ such that $h(\overline{\mu}) = h(\overline{\mu}\circ f^{-1})$. %We know that $P^{\Z}$ is identified with the base $P^{\Z} \times \{1\}$ of $\Omega$ (see \cite{MR858689})
Moreover, for a Borel subset $B$ of $ P_{\beta}^{\Z} $,
\begin{equation}
 \nu (B) = \dfrac{ \overline{\mu}( B \times \{1\} )}{\overline{\mu}(P_{\beta}^{\Z} \times \{1\} )}.\label{indm}
\end{equation}
Since $ P_{\beta} $ is positive recurrent, there is a unique measure $\overline{\mu}$ with entropy $ \log|\beta| $ on the tower $(\Omega, T) $ which induces the unique invariant probability measure $\nu $ on $P_{\beta}^{\Z}$ such that:
\begin{equation}
\nu([x])= \frac{1}{|\beta|^{l(x) }} \text{where $ x \in P_{\beta} $}.
\end{equation}
%From Proposition 2.17 of \cite{MR858689}
Note that $[x]\times\{i\}=\{((x_n)_{n\in \Z}, i) | x_n\in P_{\beta} \text{ and } x_0=x \}$. Moreover, $[x] = f^{-1}([x]\times \{1\})$. Denote by $\mu_{\beta}=\overline{\mu}\circ f^{-1}$ the measure of maximal entropy on $P_{\beta}^{\Z}$.
\begin{equation}
 \nu ([x]) \overline{ \mu } (P_{\beta}^{\Z}\times \{1\}) = \overline{ \mu }( [ x ] \times \{1\} ).
 \label{(27)}
\end{equation}
That is
\begin{equation}
 \mu_{\beta}([x]) = \nu([x]) \overline{\mu}(P_{\beta}^{\Z}\times\{1\}).
\end{equation}

Since $ \Omega = \bigcup\limits_{ x \in P_{\beta} } \left(\underset{i=1}{\overset{l( x )}{\cup}} [ x ] \times \{ i \} \right) $, it follows that:
\begin{equation}
\begin{aligned}
 1 &= \sum \limits_{ x \in P_{\beta}} \underset{i = 1}{\overset{ l( x ) }{\sum}}\overline{\mu}\left([x] \times \{i\} \right) \\
   &= \sum \limits_{ x \in P_{\beta} } \underset{i = 1}{\overset{ l( x ) }{\sum}}\overline{\mu}\left(T^{-i+1}( [x]\times\{i\}) \right) \\
   &= \sum \limits_{ x \in P_{\beta} } \underset{i = 1}{\overset{ l( x ) }{\sum}}\overline{\mu}\left([x] \times \{1\} \right) \\
   &= \sum \limits_{ x \in P_{\beta}} l( x ) \overline{\mu}\left([x] \times \{1\} \right).
\label{(28)}
\end{aligned}
\end{equation}
Moreover, the average length of $ P_{\beta} $ with respect to the measure $\nu$ is:
 \begin{equation}
\begin{aligned}
 L(P_{\beta}, \nu ) &= \sum \limits_{ x \in P_{\beta} } l( x ) \nu ([x]) \\
                   &= \dfrac{ 1 }{ \overline{ \mu }(P_{\beta}^{*}\times\{1\})}\sum\limits_{ x \in P_{\beta}} \overline{\mu}([x] \times \{1\} ) \\
                   &=\dfrac{ 1 }{ \overline{ \mu }(P_{\beta}^{*} \times \{1\} ) } \\
                   &=\sum\limits_{x\in P_{\beta}}\dfrac{l(x)}{|\beta|^{l(x)}}.
\label{(29)}                   
\end{aligned}
\end{equation}
And thus:
\begin{equation}
 \mu_{\beta} ( [ x ] ) = \left( |\beta|^{l( x )} \sum \limits_{ x \in P_{\beta} } \dfrac{ l(x)  } {|\beta|^{l(x)}} \right)^{-1}.
 \label{(30)}
\end{equation}
From Remark \ref{remMix}, $ \mu_{\beta}$ is the Champernowne measure of $P_{\beta}$.
\end{preu} 

\begin{preu}\textbf{ of Theorem \ref{Theo 4}} \\
From Corollary \ref{Support} and Corollary \ref{support 2}, a measure of maximal entropy is supported by $P_{\beta}^{\Z}$. That is, two measures $\mu_1$ and $\mu_2$ of maximal entropy coincide on $ S_{\beta} \setminus P_{\beta}^{\Z}$. Moreover, from \cite{MR858689}, since $P_{\beta}$ is a positive recurrent prefix code, there is a unique measure of maximal entropy on $P_{\beta}^{\Z}$. Thus, $\mu_1 $ and $\mu_2$ coincide on $P_{\beta}^{\Z}$ too. It follows that there is a unique measure of maximal entropy on $ S_{\beta}$.
\end{preu}

By Theorem \ref{Th tau}, for all real number $ x \in ] -\infty, -\gamma_0 ] $, there is a unique $ \beta_n $ in each interval of the type $ ] \gamma_n, \gamma_{n+1} ]$ such that $ x = \Upsilon(\beta_n) $. The system $ S_x $ is coded by $ C_x $ and the measure of maximal entropy is supported by a subset coded by $ \phi^{n+1}(C_x)$. That is, the support of the measure of maximal entropy is $ \bigcup\limits_{k \geq 0} \sigma^k(\phi^{n+1}(S_x))$ which is the interval of sequences delimited by $d_1d_2\cdots$ and $u_nd_1d_2 \cdots $. We set 
\begin{equation*}
\begin{aligned}
G_{\beta_n} &= \phi^{n+1}(S_x) \\
            &= \{ (x_i)_{i \geq 1} \in S_{\beta_n} | d^{*}(l_{\beta_n}, \beta_n) \preceq (x_i)_{i \geq 1} \preceq u_n d^{*}(l_{\beta_n}, \beta_n) \}
\end{aligned}
\end{equation*}
Since any sequence in $ G_{\beta_n}$ starts by the string $ u_n$, it follows that: 
\begin{equation*}
\begin{aligned}
\sigma^{l(u_n)} G_{\beta_n} &= \phi^{n+1}(S_x) \\
                            &= \{ (x_i)_{i \geq 1} \in S_{\beta_n} | d^{*}(l_{\beta_n}, \beta_n) \preceq (x_i)_{i \geq 1} \preceq (u_{n-1})^{2m_1}u_n (u_{n-1})^{2m_2}u_n(u_{n-1})^{2m_3} \cdots \} \\
                            &= \{ (x_i)_{i \geq 1} \in S_{\beta_n} | d^{*}(l_{\beta_n}, \beta_n) \preceq (x_i)_{i \geq 1} \preceq d(s_{l(u_n)}, \beta_n) \} \\
                            &= H_{\beta_n}.
\end{aligned}
\end{equation*}
And $ G_{\beta_n} \subset H_{\beta_n}$. Thus 
\begin{equation*}
\begin{aligned}
\bigcup\limits_{k \geq 0} \sigma^k(\phi^{n+1}(S_x)) & = \bigcup\limits_{k \geq 0} \sigma^k G_{\beta_n} \\
                                                    & = \bigcup\limits_{k \geq 0} \sigma^k H_{\beta_n}
\end{aligned}
\end{equation*} 
 The following sequences are ordered in ascending order (in the sense of the alternating order): 
\begin{align*}
& u_n(u_{n-1})^{2m_1}u_n(u_{n-1})^{2m_2}u_n(u_{m-1})^{2m_3} \cdots && = d(s_0, \beta_n) \\
& (u_{n-1})^{2m_1}u_n (u_{n-1})^{2m_2}u_n(u_{n-1})^{2m_3} \cdots   && = d(s_{l(u_n)}, \beta_n)\\
& (u_{n-1})^{2m_1-1}u_n(u_{n-1})^{2m_2}u_n(u_{n-1})^{2m_3}\cdots   && = d(s_{l(u_nu_{n-1})}, \beta_n)\\
& u_{n-2}u_{n-2}(u_{n-1})^{2m_1}u_n(u_{n-1})^{2m_2}u_n (u_{n-1})^{2m_3}\cdots && = d(s_{l(u_{n-1})}, \beta_n)\\
& u_{n-2} (u_{n-1})^{2m_1}u_n(u_{n-1})^{2m_2}u_n(u_{n-1})^{2m_3}\cdots && = d(s_{l(u_{n-1}u_{n-2})}, \beta_n) \\
& u_{n-3}u_{n-3} (u_{n-1})^{2m_1-2}u_n(u_{n-1})^{2m_2}u_n(u_{n-1})^{2m_3}\cdots &&=d(s_{l(u_nu_{n-1}u_{n-2})}, \beta_n) \\
& u_{n-3}u_{n-3}(u_{n-1})^{2m_1-1}u_n(u_{n-1})^{2m_2}u_n(u_{n-1})^{2m_3}\cdots &&= d(s_{l(u_{n}u_{n-2})},\beta_n) \\
& u_{n-3}u_{n-3}u_{n-2}u_{n-2}(u_{n-1})^{2m_1}u_n(u_{n-1})^{2m_2}u_n(u_{m-1})^{2m_3}\cdots && = d(s_{l(u_{n-2})}, \beta_n)\\
& u_{n-3}u_{n-2}u_{n-2} (u_{n-1})^{2m_1}u_n(u_{n-1})^{2m_2}u_n(u_{m-1})^{2m_3} \cdots && = d(s_{l(u_{n-2}u_{n-3})}, \beta_n) \\
& u_{n-3}(u_{n-1})^{2m_1-1}u_n(u_{n-1})^{2m_2}u_n(u_{n-1})^{2m_3}\cdots &&= d(s_{l(u_{n}u_{n-2}u_{n-3})},\beta_n)\\
& u_{n-3}(u_{n-1})^{2m_1-2}u_n(u_{n-1})^{2m_2}u_n(u_{n-1})^{2m_3}\cdots &&=d(s_{l(u_nu_{n-1}u_{n-2}u_{n-3})}, \beta_n)\\
& u_{n-4}u_{n-4} (u_{n-1})^{2m_1}u_n(u_{n-1})^{2m_2}u_n (u_{n-1})^{2m_3}\cdots && = d(s_{l(u_{n-1}u_{n-2}u_{n-3})}, \beta_n)\\
& u_{n-4}u_{n-4} u_{n-2} (u_{n-1})^{2m_1}u_n(u_{n-1})^{2m_2}u_n (u_{n-1})^{2m_3}\cdots && = d(s_{l(u_{n-1}u_{n-3})}, \beta_n)\\ 
& u_{n-4}u_{n-4}u_{n-3}u_{n-3}(u_{n-1})^{2m_1-2}u_n(u_{n-1})^{2m_2}u_n(u_{n-1})^{2m_3}\cdots  && = d(s_{l(u_nu_{n-1}u_{n-3})}, \beta_n) \\
& u_{n-4}u_{n-4}u_{n-3}u_{n-3}(u_{n-1})^{2m_1-1}u_n(u_{n-1})^{2m_2}u_n(u_{n-1})^{2m_3}\cdots && = d(s_{l(u_nu_{n-3})}, \beta_n)\\
& u_{n-4}u_{n-4}u_{n-3}u_{n-3}u_{n-2}u_{n-2} (u_{n-1})^{2m_1}u_n(u_{n-1})^{2m_2}u_n(u_{m-1})^{2m_3} \cdots && = d(s_{l(u_{n-3})}, \beta_n) \\
& \vdots && \vdots
\end{align*}
In light of the above, we remark that for a pair of nonnegative integers $ k $ and $ i $ such that $ k< n $, $ i< l(u_{k-1})$, the sequences $ \sigma^{l(u_nu_k)+i}d^{*}(l_{\beta_n}, \beta_n)$, $\sigma^{l(u_k)+i}d^{*}(l_{\beta_n}, \beta_n) $, $\sigma^{l(u_ku_{k-1})+i}d(l_{\beta_n}, \beta_n)$ and $ \sigma^{l(u_nu_ku_{k-1})+i}d(l_{\beta_n}, \beta_n)$ are the endpoints of two consecutive intervals of words in the union $ \bigcup\limits_{k \geq 0}\sigma^kH_{\beta_n}$: $ \sigma^{l(u_k)+i}H_{\beta_n}$ and  $ \sigma^{l(u_ku_{k-1})+i}H_{\beta_n}$.

We see that $ d(s_{l(u_nu_k)}+i, \beta_n)$, $ d(s_{l(u_k)+i}, \beta_n)$, $ d(s_{l(u_ku_{k-1})+i}, \beta_n)$ and $ d(s_{l(u_nu_ku_{k-1})+i}, \beta_n) $ are endpoints of consecutive intervals of $ \bigcup\limits_{k \geq 0}\sigma^kH_{\beta_n}$.

 In the right one-sided $\beta$-shift,  $ \bigcup\limits_{k \geq 0} \sigma^k(\phi^{n+1}(S^d_x))$ is the set of $\beta$-representations of $ \bigcup\limits_{k \geq 0} T_{\beta}^k ([l_{\beta}, t_n ] )$. It can be established that this union corresponds to the support of the measure of maximal entropy on the dynamical system $ (I_{\beta}, T_{\beta} ) $ determined and closely studied in \cite{MR2974214}. First, let's point out that $ \bigcup\limits_{k \geq 0} T_{\beta}^k ([l_{\beta}, t_n ] ) =  \bigcup\limits_{k \geq 0} T_{\beta}^k ([s_0, s_{l(u_n)} ] )$, where $ s_i = T_{\beta}^i(l_{\beta})$ with $ i \in \N$. Indeed, $ T_{\beta}^{l(u_n)}([s_0, t_n ]) = [s_0, s_{l(u_n)} ] $ and $ [s_0, t_n] \subset [s_0, s_{l(u_n)} ] $.

\section{Gaps on the negative beta-shift}

Let $\beta \in ]-\gamma_n, -\gamma_{n+1} ]$.  
The phenomenon of gaps on $ I_{\beta} = [l_{\beta}, r_{\beta} ) $ was closely studied in \cite{MR2974214}. In this section, we carry out an analogous study on the $\beta $-shift.

\begin{defi}
A word $ v \in L_{\beta} $ is intransitive if there exists $ u \in L_\beta$ such that for any $ w $ in $ L_{\beta}$, $ vuw \not\in L_{\beta} $. 
\end{defi}
We can see an intransitive word as a word which does not belong to the language of the support of a measure of maximal entropy.

The following result is obvious.
\begin{prop}\label{pp5}
Let $ \beta $ be a real number such that

\begin{equation}
  d^{*}(l_{\beta}, \beta) = u_n(u_{n-1})^{2k_1}u_n(u_{n-1})^{2k_2}u_n (u_{n-1})^{2k_3} \cdots.
 \end{equation}
An admissible word is intransitive if it contains one of the following sequences: 

\begin{align*}
 &\sigma^i(u_{m-2}) u_{m-1}u_m &&\text{ with $ m > 0$, $ 0 \leq i < |u_{m-2} | $,}\\
 &\sigma^i(u_{m-1}) u_{m-1}u_{m-1}u_{m-1} &&\text{ with $ m \geq 0 $, $ 0 \leq i< |u_{m-1}|$, }\\
 &\sigma^i(u_{m-1}) u_{m-1} \cdots u_{n-2} u_{n-2} (u_{n-1})^{2k_1+1} u_n &&\text{ with $ m \geq 0 $, $ 0\leq i < |u_{m-1}|$.}
 \end{align*}
with $ u_{-1} = 0 $.
\end{prop}
The words listed in Proposition \ref{pp5} are forbidden in the language of the support of a measure of maximal entropy.

It is easy to see that an admissible word $ x $ starting with $ \sigma^i(u_k) $ contains an intransitive word if and only if it is taken between
\begin{equation}
\sigma^i(u_k)u_k u_{k+1}u_{k+1}\cdots u_{n-2} u_{n-2} (u_{n-1})^{2k_1}u_n (u_{n-1})^{2k_2}\cdots = \sigma^i(d^{*}(l_{\beta}, \beta))
\end{equation}
and\begin{equation}
\sigma^i(u_k) u_{k+1}u_{k+1}\cdots u_{n-2} u_{n-2} (u_{n-1})^{2k_1}u_n (u_{n-1})^{2k_2}\cdots = \sigma^{|u_k|+i}(d^{*}(l_{\beta}, \beta)).
\end{equation}

\begin{theo}\label{InW1}
Let $ \mu $ be an ergodic measure on the symbolic system $ (X, T)$ and $ L $ its language.
Consider two words $ u $ and $ t $ of  $ L $ such that $ \forall a \in L $, $ uat \not\in L$ ( $ t$ is intransitive if there is such a word $ u $). Then, $ \mu(_0[t]) = 0 $ or $ \mu(_0[u]) = 0$.
\end{theo}

\begin{preu}
When a measure $ \mu $ is ergodic, almost every point is generic (see Proposition (5.9) of \cite{MR0457675}). Thus, if $\mu(_0[u])$ and $ \mu(_0[t]) $ are not equal to zero and if $(x_n)_{n\in\Z}$ is generic for $ \mu $, there exist infinitely many words $ u $ and $ t $ in the sequence $ (x_n)_{n \geq 1 }$ and thus a word $ a $ in $ L $ such that $ u a t \in L $ (and a word $ b $ of $ L $ such that $ tbu \in L $).
\end{preu}

We deduce the following theorem.
\begin{theo}\label{InW2}
Let $\beta<-1 $ and $ \mu $ be a measure of maximal entropy on the negative $ \beta$-shift; $ \mu (_0[x])>0 $ whenever $ x $ can be decomposed into product of words of the positive recurrent prefix (or suffix) code of its support. Let $ t $ be an intransitive word of $ L_{\beta}$. Then $\mu(_0[t]) = 0$.
\end{theo}

\section*{}

In summary, we have seen that for each case studied, there exists a unique $\sigma$-invariant measure of maximal entropy. Considering the one side $\beta$-shift, the results remain valid. If $ d(l_{\beta}, \beta) $ is periodic with odd period, $ S_{\beta}$ and $\tilde{S}_{\beta}$ have the same entropy. The negative $\beta$-shift $ S_{\beta}$ is the union of $\tilde{S}_{\beta}$ which is intrinsically ergodic and the $\sigma$-invariant subset of words ending with $ d(l_{\beta}, \beta) $. When $ \beta $ is between $ -\frac{1+\sqrt{5}}{2} $ and -1, the system is not transitive. But $S_{\beta}$ remains intrinsically ergodic. In \cite{article}, an example of subsystem of $ S_{\beta}$ not intrinsically ergodic is given by: $ X = \{1^{\infty} \} \cup \{1^n2^{\infty}: n\geq 1 \} \cup \{2^{\infty} \}$. This subshift corresponds to $\{0^{\infty} \}\cup \{ 0^n1^{\infty}: n \geq 1 \} \cup \{1^{\infty} \} $ according to our definition of the negative $\beta$-transformation. It is easy to see that this subshift is contained in all negative $\beta$-shifts. This example shows that in the intrinsically ergodic dynamical system, we can find subsystems which do not have this property. But it is necessary to attach the condition to this subsystem to have an entropy strictly smaller than the entropy of the system.

%\begin{center}
%\psset{xunit=0.01cm,yunit=1cm}
%\begin{pspicture}(-325,-0.5)(450,0.5)
 %\psline{->}(-325,0)(450,0)
 %\multido{\n=-300+100}{8}
 %{\psline[linewidth=1.2pt](\n,-0.2)(\n,0.2)
% \uput[d](\n,-0.27){$\n$}}
%\end{pspicture}
%\end{center}

%\begin{center}
% \pspicture(-5.5,-0.5)(5.5,1)
 %   \psline{->}(-6,.5)(6,.5) \multirput(-5.5,.5)(1,0){11}{+}
  %  \uput[90](-5.5,0.5){}
   % \uput[90](-4.5,0.5){$-\gamma_0$}
    %\uput[90](-3.5,0.5){$-\gamma_1$} \uput[90](-2.5,0.5){$-\gamma_2$}
    %\uput[90](-1.5,0.5){$-\gamma_3$}  \uput[90](-0.5,0.5){$-\gamma_4$}
    %\uput[90](0.5,0.5){$\cdots$}    \uput[90](1.5,0.5){$-\gamma_n$}
   % \uput[90](2.5,0.5){$-\gamma_{n+1}$}   \uput[90](3.5,0.5){$\cdots$}
  % \uput[90](4.5,0.5){-1}
 %  \psdots[dotstyle=+,dotangle=45,dotsize=0.2](-3.5,0.5)\uput[-90](-3.5,0.5){$B$}
%    \psdots[dotstyle=+,dotangle=45,dotsize=0.2](1.5,0.5)\uput[-90](1.5,0.5){$A$}
%\endpspicture
 %\end{center}

%\begin{tikzpicture}
%\node[above right] at (-1,0) {Expression};
%\node[element] (T) at (4,0) {Terme};
%\node[terminal] (p) at (4,-1) {$+$};
%\draw[fleche] (0,0) -- (T); \draw[fleche] (T) -- (8,0);
%\draw[fleche] (6,0) |- (p); \draw[fleche] (p) -| (2,0);
%\end{tikzpicture}

\bibliographystyle{unsrt}
\bibliography{Mabiblio}

\end{document}